\documentclass[11pt,a4paper,reqno]{amsart}
\usepackage[applemac]{inputenc}
\usepackage[T1]{fontenc}
\usepackage{amsmath}
\usepackage{amsthm}
\usepackage{amsfonts}
\usepackage{amssymb}
\usepackage{graphicx}
\usepackage{amsbsy}
\usepackage{mathrsfs}
\usepackage{bbm}
\newtheorem{theorem}{Theorem}[section]
\newtheorem{lemma}[theorem]{Lemma}

\newtheorem{proposition}[theorem]{Proposition}

\newtheorem{definition}[theorem]{Definition}
\theoremstyle{remark}
\newtheorem{remark}[theorem]{\it \bf{Remark}\/}

\numberwithin{equation}{section}
\catcode`@=11
\def\section{\@startsection{section}{1}%
  \z@{1.5\linespacing\@plus\linespacing}{.5\linespacing}%
  {\normalfont\bfseries\large\centering}}
\catcode`@=12
\newcommand{\be}{\begin{equation}}
\newcommand{\ee}{\end{equation}}
\newcommand{\bea}{\begin{eqnarray}}
\newcommand{\eea}{\end{eqnarray}}
\newcommand{\bee}{\begin{eqnarray*}}
\newcommand{\eee}{\end{eqnarray*}}

\def\pa{\partial}

\def\RR{\mathbb{R}}

\def\fref#1{{\rm (\ref{#1})}}

\def\G{{\Gamma}}

\catcode`@=11
\def\supess{\mathop{\operator@font Sup\,ess}}
\catcode`@=12

\def\bt{\tilde{b}}

\def\RR{\mathbb{R}}

\def\e{\varepsilon}

\def\fref#1{{\rm (\ref{#1})}}

\def\R2+{\RR ^2_+}

\def\lsl{\frac{\lambda_s}{\lambda}}

\def\pa{\partial}

\def\lim{\mathop{\rm lim}}

\def\sup{\mathop{\rm sup}}

\def\l{\lambda}

\def\log{{\rm log}}

\def\lsl{\frac{\lambda_s}{\lambda}}

\def\Psih{\hat{\Psi}}

\def\qbt{\tilde{Q}_b}

\def\tt{\tilde{T}}

\def\pa{\partial}

\def\eh{\hat{\e}}

\def\pa{\partial}

\def\Psit{\tilde{\Psi}}
\def\qbh{\hat{Q}_b}
\def\wh{\hat{w}}

\title[]{Type II blow up for the four dimensional energy critical semi linear heat equation }

\author[R. Schweyer]{R\'emi Schweyer}
\address{Institut de Math\'ematiques de Toulouse, Universit\'e Paul  Sabatier, Toulouse, France}
\email{remi.schweyer@math.univ-toulouse.fr}

\begin{document}
\maketitle
\begin{abstract}
We consider the energy critical four dimensional semi linear heat equation $\pa_tu-\Delta u-u^3=0$. We show the existence of type II finite time blow up solutions and give a sharp description of the corresponding singularity formation. These solutions concentrate a universal bubble of energy in the critical topology $$u(t,r)-\frac{1}{\l}Q\left(\frac{r}{\lambda(t)}\right)\to u^*\ \ \mbox{in}\ \ \dot{H^1}$$ where the blow up profile is given by the Talenti Aubin soliton $$Q(r)=\frac{1}{1+\frac{r^2}{8}},$$ and with speed $$\lambda(t)\sim \frac{T-t}{|\log(T-t)|^2} \ \ \mbox{as}\ \ t\to T.$$ Our approach uses a robust energy method approach developped for the study of geometrical dispersive  problems \cite{RR}, \cite{MRR}, and lies in the continuation of the study of the energy critical harmonic heat flow \cite{RS} and the energy critical four dimensional wave equation \cite{HR}.
\end{abstract}


\section{Introduction}


\subsection{Setting of the problem} We consider in this paper the energy critical semi linear heat equation
\be
\label{harmonicheatflow}
\partial_t u -\Delta u -u^3=0, \ \ (t,x)\in \Bbb R\times \Bbb R^4
\ee
which is the energy critical four dimensional version of the more general problem 
\be
\label{cbenoneoneo}
\partial_t u -\Delta u -u^p=0, \ \ (t,x)\in \Bbb R\times \Bbb R^N, \ \ p\geq 2^*-1
\ee where $$2^*=\frac{2N}{N-2}$$ is the Sobolev exponent.  There is an important literature devoted to the qualitative description of solutions to \fref{cbenoneoneo}, and we refer to \cite{MM1}, \cite{MM2} for a complete introduction to the history of the problem. For radial data, two type of blow up regimes are typically expected: type I blow up which corresponds to a self similar blow up, and type II blow up which displays excited blow up speeds. Such kind of type II blow up solutions were exhibited for the first time by Herrero and Velazquez \cite{HV} using matching asymptotic procedures for a large value of $p$, and the corresponding regime displays a polynomial type blow up speed. A major breakthrough is achieved by Matano and Merle in \cite{MM1}, \cite{MM2}, where the {\it non existence of type II} blow up is shown for $$2^*-1<p<p_c$$ where $$p_c=\left\{ \begin{array}{ll}+\infty\ \ \mbox{for}\ \ N\leq 10\\ 1+\frac{4}{N-4-2\sqrt{N-1}}\ \ \mbox{for} \ \ N\geq 11\end{array}\right.,$$ and the {\it existence} of type II blow up for $p>p_c$ is proved. More precisely, such solutions are obtained as {\it threshold} dynamics between well known type I blow up solutions and global dissipative dynamics. A complete classification of these type II regimes is then completed in \cite{Mizoguchi} where quantized blow up speeds are exhibited with polynomial rates.\\
These results leave completely open the question of existence of type II blow up in the energy critical setting. In fact, in the energy critical setting and even for the parabolic problem, the maximum principle does not seem to yield enough information to control a type II blow up. The criticality of the problem is reflected by the fact that the total dissipated energy 
\be
\label{enegyontro}
E(u)=\frac12\int_{\Bbb R^N}|\nabla u(t,x)|^2dx-\frac1{p+1}\int_{\Bbb R^N}u^{p+1}(t,x)dx
\ee
is left invariant by the scaling symmetry of the problem $$u_\lambda(t,x)=\l^{\frac{N-2}{2}}u\left(\l^2 t,\l x\right), \ \ E(u_\l(t))=E(u(\l^2t)).$$ The study of {\it critical problems} has attracted a considerable attention for the past ten years in the dispersive community, in particular the study of the mass critical nonlinear Schr\"odinger equation \cite{P}, \cite{MR1}, \cite{MR2}, \cite{MR3}, \cite{MR4}, \cite{MR5} and geometric problems like wave maps, Schr\"odinger maps and the harmonic heat flow  \cite{KST}, \cite{RodSter}, \cite{RR}, \cite{MRR}, \cite{RS}. In particular, a robust {\it energy approach} is developed in \cite{RR}, \cite{MRR} to construct type II blow up solutions in the energy critical setting. This strategy is implemented in the parabolic setting in \cite{RS} and led to the construction of {\it stable} blow up dynamics with sharp asymptotics on the singularty formation for the harmonic heat flow. Note that more type II regimes for dispersive problems were obtained in \cite{KST}, \cite{KSTwave} but rely on the construction of {\it non smooth} solutions and a procedure of backwards in time integration of the flow from the singularity which are both non suitable for parabolic problems.

\subsection{Statement of the result}

We carry out in this paper the program which was implemented in \cite{HR} to adapt the study of the geometric wave equation in \cite{RR} to the semi linear cubic four dimensional wave equation. The main difficulty is the fact that the energy \fref{enegyontro} is non definite positive, and this induces a non positive  eigenvalue in the spectrum of the linearized operator close to the Talenti-Aubin stationary solution \be
\label{ground}
Q(r)=\frac{1}{1+\frac{r^2}{8}},
\ee
which is the unique up to scaling radially symmetric solution to the stationary problem
\be
\label{soiton}\Delta Q+Q^3=0, \ \ Q(r)=\frac{1}{1+\frac{r^2}8}.
\ee
This requires building our set of initial data on a suitable codimension one set, and we similarly claim in the continuation of \cite{RS} the existence of a type II blow up dynamics for the energy critical four dimensional problem:

\begin{theorem}[Existence of type II blow up in dimension $N=4$]
\label{main}
Let $Q$ be the Talenti Aubin soliton \fref{soiton}.
 Then $\forall \alpha^*>0$, there exists a radially symmetric initial data $u_0\in H^1(\Bbb R^4)$ with
\be
E(Q)<E(u_0) < E(Q) + \alpha^*
\ee
such that the corresponding solution to the energy critical focusing parabolic equation \fref{harmonicheatflow} blows up in finite time $T = T(u_0) < \infty$ in a type II regime according to the following dynamics: there exist $u^* \in \dot H^1$ such that:
\be
\label{convustarb}
\nabla\left[u(t,x) - \frac{1}{\lambda(t)}Q\left(\frac{x}{\lambda(t)}\right)\right]\to \nabla u^*\ \ \mbox{in}\ \ L^2\ \ \mbox{as}\ \ t\to T
\ee
at the speed
\be
\label{law}
\lambda(t) = c(u_0)\left(1+o(1)\right) \frac{T-t}{|\log (T-t)|^2} \ \ \ \mbox{as} \ \ \ t \rightarrow T 
\ee
for some $c(u_0) >0$.
Moreover, there hods the regularity of the asymptotic profile:
\be
\label{regularityustar}
\Delta u^* \in L^2
\ee
\end{theorem}

{\it Comments on the result}\\

{\it 1.} In dimension four, the choice $p=2^*-1=3$ is therefore the only one for which a type II blow up occurs for radial data. We have decided to focus onto the four dimensional case for the case of simplicity but the construction we propose could be addressed in a much more general setting. Let us insist also that it does not rely on the maximum principle and may therefore  be addressed in the non radial setting and for more complicated systems.\\

{\it 2.} The blow up speed \fref{law} is also the one obtained for the energy critical harmonic heat flow in \cite{RS}. Following the heuristics developed in \cite{papierheatflow}, we conjecture the existence of a sequence of {\it quantized} blow up speeds with polynomial rates corrected by suitable logarithmic factors, and \fref{law} is the {\it fundamental} which corresponds from the proof to a codimension one in some weak sense manifold of initial data.\\

The main open problem after this work is to obtain a complete classification of type II blow up for the energy critical problem, both in the radially symmetric case and the non symmetric case.\\

{\bf Aknowledgments} This work is part of the author's PhD thesis which is supported by the ERC/ANR project SWAP. The author would like to thank his advisor Pierre Rapha\"el for his guidance and advice.\\

{\bf Notations} We introduce the differential operator $$\Lambda f=f + y\cdot\nabla f\ \ (\mbox{energy critical scaling}).$$ Given a positive number $b>0$, we let 
\be
\label{defbnot}
B_0=\frac{1}{\sqrt{b}}, \ \ B_1=\frac{|\log b|}{\sqrt{b}}.
\ee
 Given a parameter $\lambda>0$, we let $$u_\lambda(r)=\frac{1}{\lambda}u(y)\ \ \mbox{with} \ \ y=\frac{r}{\lambda}.$$  We let $\chi$ as a smooth cut off function with $$\chi(y)=\left\{\begin{array}{ll}1\ \ \mbox{for}\ \ y\leq 1,\\ 0\ \ \mbox{for} \ \ y\geq 2.\end{array}\right.$$
We shall systematically omit the measure in all radial two dimensional integrals and note: $$\int f=\int_0^{+\infty}f(r)r^3dr.$$

\section{Construction of an explicit approximate solution}

The aim of this section is to construct an approximate blow up solution of \fref{harmonicheatflow}, which is close to the ground state $Q$. This approximate solution will be the dominant part of the blow up profile inside the parabolic zone.  We adapt the strategy developed in \cite{MRR}, \cite{RR}, \cite{RS}.

\par
Let u be a stationary solution of \fref{harmonicheatflow}. Let $\lambda>0$. Then, $\frac 1{\lambda} u(\frac r {\lambda})$ is also a solution. Consider now that $\lambda$ is no more a constant, but depends of time. Thus, we obtain the following equation:
\be
\label{eq1}
\lambda(t)^2 \pa_t u\left(\frac{r}{\lambda(t)}\right) - \lambda(t)\lambda'(t) \Lambda u\left(\frac{r}{\lambda(t)}\right) - \Delta u\left(\frac{r}{\lambda(t)}\right) - u^3\left(\frac{r}{\lambda(t)}\right) = 0.
\ee 
We then define a rescaled time
\be
s = \int_0^t \frac{d\tau}{\lambda^2(\tau)}.
\ee
Remark that if $\lambda(t)$ verify the law \fref{law} defined in the theorem \ref{main}, then $s(t)$ is a bijection between $[0,T[$ and $\mathbb R^+$. We also let the rescaled variable $y(t) = \frac r {\lambda(t)}$.  The equation \fref{eq1} becomes, using the new variables:
\be
\label{eq2}
 \pa_s u - \frac{\lambda_s}{\lambda} \Lambda u - \Delta u - u^3 = 0.
\ee 
As well as the parameter $\lambda(s)$, we define a new parameter $b(s)$ such that:
\bea
\label{modul1}b &=& - \frac{\lambda_s}{\lambda}\left(1 + o(1) \right), \\
\label{modul2}b_s&=& - b^2 \left(1 + o(1) \right).
\eea
The modulation laws \fref{modul1} and \fref{modul2} will be justified thereafter.  First, in the following subsection, we consider that 
\be
\label{modulcontrainte}
b_s = -b^2 \ \ \ \mbox{and}\ \ \ b+ \frac{\lambda_s}{\lambda}=0,
\ee 
b being positive. 

\subsection{Construction of explicit approximate blow up profiles}

\begin{proposition}[Construction of the approximate profile]
\label{consprofapproch}
Let $M>0$ enough large. Then, there exists a small enough universal constant $b^*(M)$, such that the following holds. Let $b \in ]0,b^*(M)[$. Then there exists profiles $T_1$, $T_2$ and $T_3$, such that
$$Q_b(y) = Q(y) + bT_1(y)+ b^2T_2(y) + b^3T_3(y) = Q(y)+ \alpha(y)$$
generates an error
\be
\label{deferreur}
\Psi_b = - b^2 (T_1 + 2b T_2)  -\Delta Q_b - (Q_b)^3+ b \Lambda Q_b
\ee
which satisfies\\
(i) Weighted bounds: 
\be
\label{controleh2erreur}
\int_{y\leq 2B_1}|H\Psi_b|^2 \lesssim b^4|\log b|^2,
\ee
\be
\label{newbornes}
\int_{y\leq 2B_1}\frac{1}{1+y^8}|\Psi_b|^2\lesssim b^6,
\ee
\be
\label{controleh4erreur}
\int_{y\leq 2B_1}|H^2\Psi_b|^2 \lesssim \frac{b^6}{|\log b|^2}.
\ee
(ii) Flux computation: Let $\Phi_M$ be given by \fref{defdirection}, then:
\be
\label{fluxcomputationone}
\frac{(H\Psi_b,\Phi_M)}{(\Lambda Q, \Phi_M)}=-\frac{2b^2}{|\log b|}+O\left(\frac{b^2}{|\log b|^2}\right).
\ee

\begin{remark} From the proof, the profiles $(T_i)_{1\leq i\leq 3}$ display a lower order dependance in b.
\end{remark}
\end{proposition}
\begin{bf}
Proof of Proposition \ref{consprofapproch}
\end{bf}
\medskip
\par
{\bf Step1 } Computation of the error \\ \\
We expand $Q_b^3$ and formulate the error $\Psi_b$  as a polynomial expression in b:
\bee
Q_b^3 &=& Q^3 + 3bQ^2T_1 + b^2 \left( 3Q^2T_2 + 3QT_1^2\right) + b^3\left( 3Q^2T_3 + 6QT_1T_2 + T_1^3 \right)\\
 &+& R_1(T_1,T_2,T_3),
\eee
where $R_1(T_1,T_2,T_3)$ is polynomial in $(T_i)_{1\leq i\leq3}$ and contains the terms of power $(b^j)_{j \geq 4}$. 
Hence,
\bea
\label{defpsib}
\nonumber \Psi_b 
 &=&b \left( HT_1 + \Lambda Q \right) \\ 
 \nonumber &+& b^2 \left( HT_2 - T_1 + \Lambda T_1 -3Q T_1^2  \right) \\ 
 \nonumber &+& b^3\left( HT_3 - 2T_2 + \Lambda T_2 - 6QT_1T_2 - T_1^3\right) \\ 
&+& b^4\Lambda T_3 + R_1(T_1,T_2,T_3)
\eea
with 
\be
\label{defh}
H = -\Delta -3Q^2 = - \Delta - V.
\ee
Morever, 
\be
\label{defV}
V(y) = \dfrac{3}{\left(1 + \frac{y^2}{8}\right)^2},
\ee
 which yields
 \be
 \label{Gamma}
V(y)=\left\{\begin{array}{ll} 3 + O\left(y^2\right) \ \ \mbox{as} \ \ y\to 0,\\ \dfrac{192}{y^4} +O\left(\dfrac{1}{y^6}\right)\ \  \mbox{as} \ \ y\to +\infty,\end{array}\right . 
\ee
and
\be
\label{comportementv}
\Lambda V = \frac{-192(3y^2-8)}{(y^2+8)^3}.
\ee
\\
\par
\begin{bf}
Step 2
\end{bf}
Construction of $T_1$ \\

The spectral structure of the Schr\"odinger operator $H$ is well known: it has a well localized non positive eigenvalue $$H\psi=-\zeta \psi, \ \ \zeta>0,$$ and a resonance at the origin induced by the energy critical scaling symmetry: $$H\Lambda Q=0, \ \ \Lambda Q\notin L^2(\Bbb R^4).$$ Hence the Green's functions of $H$ are explicit and the other solution to $H\Gamma=0$ for $y>0$ is given by:
$$ \Gamma(y)= - \Lambda Q \int_1^y\frac{dx}{x^3(\Lambda Q(x))^2} = \dfrac{y^2-8}{(y^2+8)^2}\left( \frac{y^2}{16} +6 \log y -\frac{583}{112} - \frac{4}{y^2} \right) - \dfrac{64}{ \left( y^2+8\right)^2},
 $$
 which yields
 \be
 \label{Gamma}
\Gamma(y)=\left\{\begin{array}{ll} O\left(\dfrac{1}{y^2}\right) \ \ \mbox{as} \ \ y\to 0,\\ \dfrac{1}{16}+O\left(\dfrac{\log y}{y^2}\right)\ \  \mbox{as} \ \ y\to +\infty.\end{array}\right . 
\ee
We may thus invert $H$ explicitely and the smooth solutions at the origin of $$Hu=f$$ are given by 
\be
\label{ezaripou}
u=\G(y)\int_0^y f\Lambda Q -\Lambda Q(y)\int_{0}^y f\G+c\Lambda Q(y), \ \ c\in \Bbb R.
\ee
We let $T_1$ be the solution of \be HT_1 + \Lambda Q =0,\ee
canceled in zero, which means that we choose the constant $c=0$ in \fref{ezaripou}.
 There holds the behaviors at $r\to +\infty$
\be
\label{asympT1}
\Lambda^iT_1(y) = -4\left( \log y - \frac{1}{2}+i\right) + O\left(\frac{(\log y)^2}{y^2}\right), \ \ \mbox{for} \ \ 0 \leq i \leq 2.
\ee
There holds the behaviors at $r\to 0$
\be
\label{asympT10}
\Lambda^i T_1 = O \left( y^2\right), \ \ \mbox{for} \ \ 0 \leq i \leq 2.
\ee
Hence, for $0 \leq i \leq 2$
\bee
\| \Lambda ^i T_1 \|_{L^{\infty}_{y \leq 2B_1}} &\lesssim& \log b, \\
H\Lambda ^iT_1 &\sim& \frac{8}{y^2}, \ \ \mbox{when} \ {y \to \infty}.
\eee
\par
\begin{bf}
Step 3
\end{bf}
Construction of the radiation $\Sigma_b$\\ \par
First, we can notice that the choice $b_s = -b^2$ allows to cancel the $\log y$ growth of the expression $- T_1 + \Lambda T_1 -3Q T_1^2$. We now construct a radiation term according two specifications. First, it must compensate the 1-growth of the last expression. Moreover the error induced by this term inside the parabolic zone $y\lesssim B_0$ must be sufficiently small in order not to perturb the dynamics of the blow-up. \\ \\
Let
\be
\label{cb}
c_b = \frac{64}{\int \chi_{\frac{B_0}{4}} (\Lambda Q)^2} =\frac{ 2}{ |\log b|}\left( 1 + O\left( \frac{1}{|\log b|}\right)\right)
\ee
and
\be
\label{db}
d_b = c_b\int_0^{B_0} \chi_{\frac{B_0}{4}} \Lambda Q \G = O\left( \frac{1}{b|\log b|}\right).
\ee
Let $\Sigma_b$ be the solution to 
\be
\label{defsigmab}
H\Sigma_b=c_b\chi_{\frac{B_0}{4}}\Lambda Q + d_b H[(1 - \chi_{3B_0}) \Lambda Q]
\ee
given by $$\Sigma_b(y)=\G(y)\int_0^y c_b\chi_{\frac{B_0}{4}}(\Lambda Q)^2 -\Lambda Q(y)\int_{0}^y c_{b}\chi_{\frac{B_0}{4}}\G\Lambda Q + d_b (1 - \chi_{3B_0}) \Lambda Q.$$
The choice of the constant $c_b$ and $d_b$ yield:
\bea
\label{sigmab0}
\Sigma_b = \left\{\begin{array}{ll} c_bT_1\ \ \mbox{for} \ \ y \leq \frac{B_0}{4}\\ \\
64\G \ \ \mbox{for} \ \ y \geq 6B_0
\end{array} \right. \ \ \
\Lambda \Sigma_b = \left\{\begin{array}{ll} c_b\Lambda T_1\ \ \mbox{for} \ \ y \leq \frac{B_0}{4}\\ \\
64\Lambda \G \ \ \mbox{for} \ \ y \geq 6B_0.
\end{array} \right.
\eea
Then the estimate for $\Sigma_b$ and $\Lambda \Sigma_b$ for $6B_0 \leq y \leq 2B_1$, 
\be
\label{sigmabinfty}
\Sigma_b(y) = 4 + O\left(\frac{\log y}{y^2}\right) \ \ \ \Lambda \Sigma_b(y) = 4 + O\left(\frac{\log y}{y^2}\right),
\ee
which fits the first criterion, that we fixed to construct the radiation. For $\frac{B_0}{4}\leq y \leq 6B_0$, we have:
\bea
\label{sigmabmilieu}
\nonumber \Sigma_b(y) & = &c_b\left(\frac{1}{16}+O\left(\frac{\log y}{y^2}\right)\right)\left[\int_0^y \chi_{\frac{B_0}{4}}(\Lambda Q)^2 \right]-c_b\Lambda Q(y)\int_1^{y}O(x)dx\\
& = & 4\frac{\int_0^y\chi_{\frac{B_0}{4}}(\Lambda Q)^2 }{\int \chi_{\frac{B_0}{4}}(\Lambda Q)^2}+O\left(\frac{1}{|\log b|}\right).
\eea
\bea
\label{lsigmabmilieu}
\nonumber \Lambda \Sigma_b(y) &=& \Lambda \G(y)\int_0^y c_b\chi_{\frac{B_0}{4}}(\Lambda Q)^2 -\Lambda^2 Q(y)\int_{0}^y c_{b}\chi_{\frac{B_0}{4}}\G\Lambda Q \\
\nonumber & = &c_b\left(\frac{1}{16}+O\left(\frac{\log y}{y^2}\right)\right)\left[\int_0^y \chi_{\frac{B_0}{4}}(\Lambda Q)^2 \right]-c_b\Lambda^2 Q(y)\int_1^{y}O(x)dx\\
& = & 4\frac{\int_0^y\chi_{\frac{B_0}{4}}(\Lambda Q)^2 }{\int \chi_{\frac{B_0}{4}}(\Lambda Q)^2}+O\left(\frac{1}{|\log b|}\right).
\eea
Similarly ,
\bee
\label{lsigmabmilieu}
 \Lambda^2 \Sigma_b(y) = 4\frac{\int_0^y\chi_{\frac{B_0}{4}}(\Lambda Q)^2 }{\int \chi_{\frac{B_0}{4}}(\Lambda Q)^2}+O\left(\frac{1}{|\log b|}\right).
\eee
The equation \fref{defsigmab} and the cancellation $H\Lambda Q=0$ yield the bounds:
\be
\label{boundhigmab}
\int|H\Sigma_b|^2\lesssim \frac{1}{|\log b|},
\ \ \
\int\frac{1}{1+y^8}|\Sigma_b|^2 \lesssim b^2,
\ \ \
 \int|H^2\Sigma_b|^2\lesssim \frac{b^2}{|\log b|^2}.
\ee
We will see that this bounds respect the second criterion for the conception of the radiation. Furthermore we will see the importance of the term $c_b$, which modify the modulation equation of $b$. It becomes:
\be
\label{wxwxwx}
b_s = -b^2\left(1+ \frac{2}{|\log b|}\right).
\ee
After reintegration, this equation gives the expected blow-up speed \fref{law}.
\\
\par
\begin{bf}
Step 4
\end{bf}
Construction of $T_2$\\ \\
Define
\bea
\label{defsigma2}
\Sigma_2 &=& \Sigma_b + T_1 - \Lambda T_1 - 3QT_1^2.
\eea
The profile $T_2$ will be defined later as the suitable output of H for the argument $\Sigma_2$. Estimate $\Sigma_2$ before the choice of $T_2$. For $y \leq 1$,
\be
\Sigma_2 \lesssim y^2.
\ee
For $1 \leq y \leq 6B_0$
\bea
\Sigma_2 \nonumber &=& 4\left(\frac{\int_0^y\chi_{\frac{B_0}{4}}(\Lambda Q)^2 }{\int \chi_{\frac{B_0}{4}}(\Lambda Q)^2} - 1\right) + O \left( \frac{(\log y)^2 }{y^2}\right)+O\left(\frac{1}{|\log b|}\right)\\
&\lesssim&  \frac{1 + \log (y\sqrt b)}{|\log b|}.
\eea
According to the choice of the modulation law for $b$ and the conception of the radiation, we have, for $y \geq 6B_0$
\be
\Sigma_2 \lesssim \frac{(\log y)^2}{y^2}.
\ee
Hence, we obtain:
\bea
\label{sigma2bord}
|\Sigma_2| \lesssim \frac{y^2}{1+y^2}\left({\bf 1}_{ y \leq 1} + \frac{1 + \log (y\sqrt b)}{|\log b|} {\bf 1}_{1 \leq y \leq 6B_0}\right) + \dfrac{(\log y )^2}{y^2} {\bf 1}_{y \geq 6B_0}.
\eea
We have the same bound for $\Lambda \Sigma_2$ and for $\Lambda^2 \Sigma_2$. We now let $T_{2}$ be the solution to 
\be
\label{eqtzerotwo}
HT_{2}=\Sigma_{2}
\ee given by 
\be
\label{formulatzerotwo}
 T_{2}(y)=\G(y)\int_0^y\Sigma_{2}\Lambda Q - \Lambda Q(y) \int_{0}^y\Sigma_{2}\G.
\ee
We derive from \fref{sigma2bord} the bounds: 
\bea
\label{decayttwo}
\forall y\leq 2B_1, \ \ |\Lambda^iT_{2}(y)|&\lesssim&  \frac{y^4}{1+y^4}\left({\bf 1}_{ y \leq 1} + \frac1{b|\log b|} {\bf 1}_{y \geq 1}\right), \ \ 0 \leq i \leq 1\\
\label{roughT2}|T_2(y)| &\lesssim& y^2.
\eea
With an explicit calculus, we prove that for any function $f$: 
\be
\label{calculrealcinr}
H\Lambda f=2Hf+\Lambda Hf-\Lambda Vf.
\ee
Hence,
\bea
\label{defhlT2}
H(\Lambda T_2) =  2\Sigma_2 + \Lambda \Sigma_2 -\Lambda VT_2
\eea
and
\bea
|H(\Lambda T_2)| \lesssim \frac{y^2}{1+y^2}\left({\bf 1}_{ y \leq 1} + \frac{1 + \log (y\sqrt b)}{|\log b|} {\bf 1}_{1 \leq y \leq 6B_0}\right) + \dfrac{(\log y )^2}{y^2} {\bf 1}_{y \geq 6B_0}.
\eea
\par
\begin{bf}
Step 5
\end{bf}
Construction of $T_3$\\
\\
In the same way as before, we define 
\be
\label{defsigma3}
\Sigma_3 = - 2T_2 + \Lambda T_2 -6QT_1T_2 -T_1^3.
\ee
Notice that we haven't to conceive a second radiation term.We estimate from \fref{asympT1} and \fref{decayttwo}
\be
\label{sigma3}
\forall y\leq 2B_1, \ \ |\Sigma_3(y)|\lesssim  \frac{y^4}{1+y^4}\left({\bf 1}_{ y \leq 1} + \frac1{b|\log b|} {\bf 1}_{y \geq 1}\right)
\ee
and
\be
\label{lsigma3}
\forall y\leq 2B_1, \ \ |\Lambda \Sigma_3(y)|\lesssim  \frac{y^4}{1+y^4}\left({\bf 1}_{ y \leq 1} + \frac1{b|\log b|} {\bf 1}_{y \geq 1}\right).
\ee
We then let $T_3$ be the solution to
 \be
 \label{eqtzerothree}
 HT_{3}=\Sigma_{3}
 \ee given by:
 \be
 \label{defT3}
 T_{3}(y)=\G(y)\int_0^y\Sigma_{3}\Lambda Q-\Lambda Q(y)\int_{0}^y\Sigma_{3}\G.
\ee
Hence,
\be
\label{defsigmaT3}
\Lambda T_{3}(y)=\Lambda \G(y)\int_0^y\Sigma_{3}\Lambda Q-\Lambda^2 Q(y)\int_{0}^y\Sigma_{3}\G.
\ee
We estimate from \fref{sigma3}
\bea
\label{T3}
\forall y\leq 2B_1, \ \ |\Lambda^iT_3(y)|&\lesssim & \frac{y^6}{1+y^4}\left({\bf 1}_{ y \leq 1} + \frac1{b|\log b|} {\bf 1}_{y \geq 1}\right), \ \ 0 \leq i \leq 1\\
\label{roughT3}
 |T_3(y)|& \lesssim& y^2(1+y^2).
\eea
Finally with \fref{calculrealcinr},
\bea
\label{defhlT3}
H(\Lambda T_3) =  2\Sigma_3 + \Lambda \Sigma_3 - \Lambda V T_3
\eea
and
\be
\label{hlt3}
\forall y\leq 2B_1, \ \ |H\Lambda T_3(y)|\lesssim  \frac{y^4}{1+y^4}\left({\bf 1}_{ y \leq 1} + \frac1{b|\log b|} {\bf 1}_{y \geq 1}\right).
\ee
We have thus the bounds for $i=0,1$, using \fref{hlt3}, \fref{T3}, \fref{sigma3} and \fref{lsigma3}:
\be
\label{estththreone}
\int_{y\leq 2B_1}|H\Lambda^i T_3|^2\lesssim \int_{y\leq 2B_1}\frac{1}{b^2|\log b|^2}\lesssim \frac{B_1^4}{b^2|\log b|^2}\lesssim \frac{|\log b|^2}{b^4}, 
\ee
\be
\label{lT3y8}
\int_{y\leq 2B_1}\frac{1}{1+y^8}|\Lambda^i T_3|^2 \lesssim \frac{1}{b^2| \log b|^2} \int_{y \leq 2B_1} \frac{1}{1+y^4}\lesssim \frac{1}{b^2}.
\ee
The crucial bounds to control the error at $H^4$ level is:
\be
\label{estththreonebis}
\int_{y\leq 2B_1}|H^2\Lambda^i T_3|^2\lesssim \frac{1}{b^2|\log b|^2}\ \ \ \mbox{for}\ \ \ \ 0 \leq i \leq 1.
\ee
We now prove this bound. A rough estimate looses the huge gain $\frac{1}{|\log b|^2}$, and we need to be more precise. From \fref{defhlT3},
 \be
 \label{eqhtrop}
 H^2(\Lambda T_3)=H(2\Sigma_3+\Lambda \Sigma_3)+O\left(\frac{1}{1+y^2}\right).
 \ee
We use again \fref{calculrealcinr}, together with \fref{defsigma3}, \fref{eqtzerotwo}, to obtain:
$$H\Sigma_3=-2HT_2+H\Lambda T_2+O\left(\frac{|\log y|^3}{y^2}\right)=\Lambda \Sigma_2+O\left(\frac{|\log y|^5}{y^2}\right),$$ $$H\Lambda \Sigma_3=2\Lambda \Sigma_2+\Lambda^2\Sigma_2+O\left(\frac{|\log y|^5}{y^2}\right)$$ 
and injecting this into \fref{eqhtrop} with \fref{sigma2bord} yields:
 \bee
 \int_{y\leq 2B_1}|H^2(\Lambda T_3)|^2 & \lesssim &\int_{y\leq 2B_1} \left|\frac{y^2}{1+y^2}\left({\bf 1}_{ y \leq 1} + \frac{1 + \log (y\sqrt b)}{|\log b|} {\bf 1}_{1 \leq y \leq 6B_0}\right) + \dfrac{(\log y )^2}{y^2} {\bf 1}_{y \geq 6B_0}\right|^2\\
 & \lesssim & \frac{1}{b^2|\log b|^2}
 \eee
and \fref{estththreonebis} is proved.\\
\par
\begin{bf}
Step 6
\end{bf}
Estimate the error\\
\\
We are in position to estimate the error $\Psi_b$. According to our construction, we have from \fref{defpsib}
\be
\label{newpsib}
\Psi_b = b^2 \Sigma_b + b^4 \Lambda T_3 + R_1(T_1,T_2,T_3)
\ee
We then study the last term, the others being already estimated.
The bounds \fref{asympT10}, \fref{asympT1}, \fref{roughT2} and \fref{roughT3} yield the bound for $y \leq 2B_1$, $0\leq i \leq 4$, and $0\leq j \leq 5$:
\bee
 \left| \frac{d^iR_1(y)}{dy^i}\right| &\lesssim&b^4\left( y^{4-i} {\bf 1}_{y \leq 1} + b^jy^{2(j+1)-i}\left(1 + |\log y|^2\right) {\bf 1}_{y \geq 1} \right) \\
&\lesssim& b^4\left( y^{4-i} {\bf 1}_{y \leq 1} + y^{2-i}|\log b|^C{\bf 1}_{y \geq 1} \right).
\eee
Hence:
\bee
\int_{y\leq 2B_1}|HR_1|^2 &\lesssim& b^8|\log b|^C\int_{y\leq 2B_1}1\lesssim b^6|\log b|^C, \\
\int_{y\leq 2B_1}\frac{1+|\log y|^2}{1+y^4}|HR_1|^2&\lesssim& b^8|\log b|^C\int_{y\leq 2B_1}\frac{1+ |\log y|^2}{1+y^4}\lesssim b^8|\log b|^C, \\
\int_{y\leq 2B_1}|H^2R_1|^2 + \frac{|R_1|^2}{1+y^8}&\lesssim& b^8|\log b|^C\int_{y\leq 2B_1}\frac{1}{1+y^4}\lesssim b^8|\log b|^C.
\eee
Injecting these bounds together with \fref{boundhigmab}, \fref{estththreone}, \fref{lT3y8}, \fref{estththreonebis} into \fref{newpsib} yields \fref{controleh2erreur}, \fref{newbornes}, \fref{controleh4erreur}.\\
We now prove the flux computation \fref{fluxcomputationone}, which will be helpful for the improved modulation equations.
\bee
\frac{(H\Psi_b,\Phi_M)}{(\Lambda Q, \Phi_M)}&=& \frac{1}{(\Lambda Q, \Phi_M)}\left[ \left(-b^2c_b\chi_{\frac{B_0}4}\Lambda Q,\Phi_M\right) + O\left(C(M)b^3\right)\right]\\
&=&-c_bb^2+O\left(C(M)b^3\right)=-\frac{2b^2}{|\log b|}+O\left(\frac{b^2}{|\log b|^2}\right).
\eee
Here we recall that $M$ large enough being chosen, we assume $|b|<b^*(M)$ so that the above claim make sense. This concludes the proof of Proposition \ref{consprofapproch}.
\subsection{Localization of the profile}

Taking a careful look at the profiles $\left(T_i \right)_{1\leq i \leq3}$, we can notice that for $y \gg B_1$, $Q$ is negligible compared to $bT_1 + b^2T_2 + b^3T_3$. Obviously, this doesn't make sense, because we look for a solution close to Q. So, we must localize the profiles, with cut-off smooth functions. For technical reasons, we use two localizations: one at $B_1$, another one at $B_0$. 
\begin{proposition}[Localization of the profile near $B_1$]
\label{localisation}
Let a $\mathcal C^1$ map $ s \mapsto b(s)$ defined on $[0,s_0]$ with a priori bound $\forall s \in [0,s_0]$,
\be
\label{bornes}
0<b(s)<b^*(M),\ \ |b_s| \leq 10 b^2.
\ee
Let the localized profile
$$\tilde{Q}_b(s,y) = Q + b \tilde{T}_1 + b^2\tilde{T}_2 + b^3\tilde{T}_3 = Q + \tilde{\alpha}
$$
where
$$\tilde{T}_i= \chi_{B_1}T_i, \ \ 1\leq i\leq 3.
$$
Then
\be
\label{eqerreurloc}
\partial_s \qbt -\Delta \qbt-\frac {\lambda_s}{\lambda} \Lambda\qbt - \qbt^3 =\mbox{Mod(t)}+ \tilde{\Psi}_b
\ee
with 
\be
\label{defmode}
Mod(t)=- \left( \frac{\lambda_s}{\lambda} +b\right) \Lambda \qbt + (b_s + b^2)(\tilde{T_1} +2b\tilde{T_2})
\ee
and where $\tilde{\Psi}_b$ satisfies the bounds on $[0,s_0]$:\\
(i) Weighted bounds: 
\be
\label{controleh2erreurtilde}
\int|H\tilde{\Psi}_b|^2 \lesssim b^4|\log b|^2,
\ee
\be
\label{newbornestilde}
\int \frac{1}{1+y^8}|\tilde{\Psi}_b|^2\lesssim b^6,
\ee
\be
\label{controleh4erreurtilde}
\int|H^2\tilde{\Psi}_b|^2 \lesssim \frac{b^6}{|\log b|^2}.
\ee
(ii) Flux computation: Let $\Phi_M$ be given by \fref{defdirection}, then:
\be
\label{fluxcomputationonebis}
\frac{(H\Psit_b,\Phi_M)}{(\Lambda Q, \Phi_M)}=-\frac{2b^2}{|\log b|}+O\left(\frac{b^2}{|\log b|^2}\right).
\ee
\end{proposition}

We introduce a second localization at $B_0$ which will relevant for $H^2$ control, see the proof of Proposition \ref{ttpptt}.

\begin{proposition}[Second localization]
\label{localisationbis}
Let a $\mathcal C^1$ map $ s \mapsto b(s)$ defined on $[0,s_0]$ with a priori bound \fref{bornes}. Let the localized profile
\be
\label{defhatqb}
\hat{Q}_b(s,y) = Q + b \hat{T}_1 + b^2\hat{T}_2 + b^3\hat{T}_3 = Q + \hat{\alpha}
\ee
where
$$\hat{T}_i= \chi_{B_0}T_i, \ \ 1\leq i\leq 3.
$$
Let the radiation:
\be
\label{defraidiation}
\zeta_b=\tilde{\alpha}-\hat{\alpha}
\ee
and the error
$$
\partial_s \qbh -\Delta \qbh-\frac {\lambda_s}{\lambda} \Lambda\qbh - \qbh^3=\widehat{Mod}(t)+ \hat{\Psi}_b
$$
with 
\be
\label{defmodebis}
\widehat{Mod}(t)=- \left( \frac{\lambda_s}{\lambda} +b\right) \Lambda \qbh + (b_s + b^2)(\hat{T_1} +2b\hat{T_2}).
\ee
Then there holds the bounds:
\be
\label{radiationnormeinfty}
\left\| \pa^i_y\zeta_b \right\|^2_{L^{\infty}} \lesssim b^{2+i} |\log b|^C,
\ee
\be
\label{radiation}
\int|H\zeta_b|^2\lesssim b^2|\log b|^2,\ \ \ \Sigma_{i=0}^2\int\frac{|\pa_y^i\zeta_b|^2}{1+y^{2(3-i)}}\lesssim b^3|\log b|^C,
\ee
\be
\label{radiationbis}
\int|H^2\zeta_b|^2 + \Sigma_{i=0}^2\int\frac{|\pa_y^i\zeta_b|^2}{1+y^{8-2i}}\lesssim b^4|\log b|^C,
\ee
\be
\label{controleh2erreurtildehat}
Supp(H\hat{\Psi}_b)\subset [0,2B_0] \ \ \mbox{and} \ \ \int|H\hat{\Psi}_b|^2 \lesssim b^4|\log b|^2.
\ee
\end{proposition}

The proof follows similar lines as in \cite{RS} and is displayed for the reader's convenience in Appendix C.

 \section{Presentation of possible solution of Theorem \ref{main}}
 \subsection{Uniqueness of the decomposition}
 We now look for a solution of \fref{harmonicheatflow} $u$, which we will decompose in the form of:
 \be
 \label{decompo}
 u = \left(\tilde Q_{b(t)}+ \e(t)\right)_{\lambda(t)} .
 \ee
 Naturally, we must fix constrains to obtain the uniqueness of this decomposition. Moreover it's crucial that the radiation term $\e$ doesn't perturb the modulation equation \fref{wxwxwx} found during the construction. We will see in the subsection devoted to the modulation equations that it's the case if we have the following inequality:
 \bee
 \int |H^2(\e(t))|^2 \lesssim \frac{b^4(t)}{|\log b(t)|^2}.
 \eee
 To control sharply the radiation term $\e$, it appears then that this one ought to be orthogonal to the kernel of $H^2$. The smooth solutions of $H^2f = 0$ are situated in $Span(\Lambda Q, T_1)$. But neither $\Lambda Q$ nor $T_1$ are in $L^2_{rad}(\mathbb R^4)$. Therefore, we use an approximation of the kernel, localizing both directions, with the smooth cut-off function $\chi_M$, where  $M>0$ is an enough large constant.  
 More precisely, we let the direction: 
 \be
 \label{defdirection}
 \Phi_M = \chi_M\Lambda Q - c_M H(\chi_M \Lambda Q)
 \ee
 with
  \bee
  c_M = \dfrac{\left(\chi_M \Lambda Q,T_1\right)}{\left(H\left( \chi_M \Lambda Q \right), T_1\right)} = c_{\chi} \frac{M^2}{4}(1+o_{M \to + \infty}(1)).
 \eee
 The second term, which is a corrective term, makes the orthogonality between $\Phi_M$ and $T_1$, and the orthogonality between $\Lambda Q$ and $H \Phi_M$, because H is a self-adjoint operator. Furthermore, in accordance with conception of this direction, we have:
 \be
 \label{intM}
 \int |\Phi_M|^2 \lesssim |\log M|,
 \ee
and the scalar products
\be
\label{prodscalar}
(\Lambda Q, \Phi_M) = (-HT_1,\Phi_M) = (\chi_M\Lambda Q,\Lambda Q) = 64 \log M (1+o_{M \to + \infty}(1)).
\ee
In the appendix, we argue that  we have coercive estimates for the operators $H$ and $H^2$ under additional orthogonality conditions. As a consequence, we fix for the radiation term $\e$ the orthogonality conditions:
\be
\label{ortho}
(\e(t),\Phi_M)=(\e(t),H\Phi_M)=0.
\ee
From a standard argument based on the implicit function theorem, these constrains give us the existence and the uniqueness of the decomposition \fref{decompo}. First, we have:
\bee
(b,\lambda) \rightarrow \left(u,\Phi_M\right) = \left(\left(\tilde Q_{b(t)}\right)_{\lambda(t)},\Phi_M\right)
\eee
is a $\mathcal C^1$ map and thus:
\bea
\nonumber
\left|\begin{array}{ll} (\frac{\pa}{\pa\l}(\qbt)_\lambda,\Phi_M)& (\frac{\pa}{\pa b}(\qbt)_\lambda,\Phi_M)\\  (\frac{\pa}{\pa\l}(\qbt)_\lambda,H\Phi_M)& (\frac{\pa}{\pa b}(\qbt)_\lambda,H\Phi_M)\end{array}\right |_{\l=1,b=0}& = & \left|\begin{array}{ll} (-\Lambda Q,\Phi_M)& 0 \\ 0 & (T_1,H\Phi_M)\end{array}\right |\\
\label{resjacobian}& = & -(\Lambda Q,\Phi_M)^2\neq 0.
\eea 
We used the orthogonality conditions mentioned at the moment of the conception of $\Phi_M$ and the following equality:
 $$\left(\frac{\pa}{\pa\l}(\qbt)_\lambda,\ \frac{\pa}{\pa b}(\qbt)_\lambda\right)|_{\l=1,b=0}=-(\Lambda Q,T_1).$$
 As long as the solution remains in a fixed small neighbourhood of $Q$ for the norm $\dot H^1$ what will be ensured for a suitable set of initial data, the implicit function theorem ensures the existence and uniqueness of the decomposition \fref{decompo}.
  \subsection{Partial differential equation verified by the radiation and suitable energies}
  \label{eed}
From now on, we always use the last decomposition. Moreover depending on whether we have the need of original variables, or rescaled variables, we shall notice the radiation term namely:
\be
 \label{decompobis}
 u = \frac 1{\lambda(t)}\left(\tilde Q_{b}+ \e \right)\left(t,\frac r{\lambda(t)}\right) = \frac 1{\lambda(t)} \tilde Q_{b(t)}\left(t,\frac r{\lambda(t)}\right) + w(t,r).
\ee
We use the occasion to recall the correspondance between both systems of variables
\bee
s(t) = \int_0^t \frac{d\tau}{\lambda^2(\tau)} \ \ \ \mbox{and} \ \ \ y = \frac{r}{\lambda(t)}.
\eee
We give also the rescaling formulas
  \bea
  \nonumber u(t,r) = \frac{1}{\lambda}v(s,y) \mbox{,}  \ \ \ \ \ \partial_t u = \frac{1}{\lambda^2}\left(\partial_s v - \frac{\lambda_s}{\lambda}\Lambda v\right)_{\lambda}.
  \eea  
  We then can inject the decomposition \fref{decompobis} with the rescaled variables in the equation \fref{harmonicheatflow} using the one of $\tilde Q_b$ \fref{eqerreurloc} and obtain the following:
  \be
  \label{eqofepsilon}
  \partial_s \varepsilon  - \frac{\lambda_s}{\lambda} \Lambda \varepsilon + H \varepsilon = F - Mod = \mathcal F,
  \ee
where we remind that
\bee 
H = - \Delta - V \mbox{,}  \ \ \ \ \ Mod(t) = - \left( \frac{\lambda_s}{\lambda} +b\right) \Lambda \qbt + (b_s + b^2)(\tilde{T_1} +2b\tilde{T_2})
 \eee
 and where we noticed
  \be
  \label{defF}
  F = -\Psit_b +L(\varepsilon) + N(\varepsilon),
  \ee
where $L$ is a linear operator coming from the difference between $H$ and $H_{B_1}$: 
 \be
\label{defLepsilon}
L(\varepsilon) = H\e - H_{B_1} \e = 3(\tilde Q_b^2 - Q^2) \varepsilon
\ee
with
\bee
H_{B_1} = - \Delta - 3 \tilde Q_b^2
\eee
  and a last purely nonlinear term:
\be
\label{defNepsilon}
N(\varepsilon) = 3\tilde Q_b \varepsilon^2 + \varepsilon ^3.
\ee
It's important to remark that we used here the localization of the profiles near $B_1$. At the end of this subsection, we will introduce in the same way some new operators with the second localization near $B_0$. Before rewriting \fref{eqofepsilon} with the original variables, introduce the suitable norms for our study:
\begin{itemize}
\item{Energy bound}
\be
\mathcal E_1 = \int|\nabla \e|^2
\ee
\item{Higher order Sobolev norms}
\be
\label{defnorme}
\mathcal E_2 = \int |H\e|^2 = \int |\e_2|^2\mbox{,}  \ \ \ \ \  \mathcal E_4 = \int |H^2\e|^2= \int |\e_4|^2
\ee
with $\e_i = H^i\e$ for $i \in \{2;4\}$.
\end{itemize} 
To work with the original variables, we would recall that we have the following notation:
\bee
f_{\lambda}(y) = \frac 1{\lambda} f(\frac r {\lambda}).
\eee
Furthermore we must adapt this notation for the potential term namely because of its quadratic nature:
\bee
\tilde V(y) = \frac 1 {\lambda^2} V(\frac r {\lambda}).
\eee
Then \fref{eqofepsilon} becomes:
  \be
  \label{eqenwini}
  \partial_t w + H_{\lambda}w = \frac 1{\lambda^2} \mathcal F_{\lambda}.
  \ee
  We define the same both functions: $w_i = H_{\lambda}^i w$ for $i \in \{2;4\}$, which verify respectively:
  \be
\label{eqw2}
\pa_t w_2 + H_{\lambda} w_2 = - \pa_t \tilde{V} w +  H_{\lambda}\left(\frac1{\lambda^2}\mathcal{F}_{\lambda}\right), 
\ee
\be
\label{eqw4}
\pa_t w_4 + H_{\lambda} w_4 = - \pa_t \tilde{V} w_2 - H_{\lambda}\left(\pa_t \tilde{V}w\right)  +  H^2_{\lambda}\left(\frac1{\lambda^2}\mathcal{F}_{\lambda}\right).
\ee
We have also by substitution:
\be
\lambda^2 \mathcal E_2 = \int |Hw|^2 = \int |w_2|^2\mbox{,}  \ \ \ \ \  \lambda^6 \mathcal E_4 = \int |H^2w|^2= \int |w_4|^2.
\ee
\\
Before closing this part, we are now getting interested us in the localization near $B_0$. Using the definition of the radiation $\zeta_b$ \fref{defraidiation}, we obtain the new unique decomposition:
\be
u = (\hat Q_b + \hat{\e})_{\lambda} \ \ \mbox{ie} \ \ \ \hat{\e} = \e + \zeta_b.
\ee  
Thus, with this localization, we have:
\be
  \label{eqofepsilonhat}
  \partial_s \eh  - \frac{\lambda_s}{\lambda} \Lambda \eh + \hat H \eh = \hat F - \widehat{Mod} = \hat{\mathcal F},
  \ee
where we remind that
\bee 
\hat H = - \Delta - \hat V \mbox{,}  \ \ \ \ \ \widehat{Mod}(t) = - \left( \frac{\lambda_s}{\lambda} +b\right) \Lambda \hat Q_{b} + (b_s + b^2)(\hat{T_1} +2b\hat{T_2})
 \eee
 and where we have noticed
  \be
  \label{defFhat}
  \hat F = -\hat \Psi_b +\hat L(\eh) + \hat N(\eh),
  \ee
where $\hat L$ is a linear operator coming from the difference between $\hat H$ and $\hat H_{B_1}$: 
 \be
\label{defLepsilon}
\hat L(\varepsilon) = \hat H\eh - \hat H_{B_1} \eh = 3(\hat Q_b^2 - \hat Q^2) \eh
\ee
with
\bee
\hat H_{B_1} = - \Delta - 3  \hat Q_b^2
\eee
  and a last purely nonlinear term:
\be
\label{defNepsilonhat}
\hat N(\eh) = 3\hat Q_b \eh^2 + \eh ^3.
\ee
We define likewise the operators $\eh_2$, $\hat w$ and  $ \hat w_2$, which come of course respectively from $\e_2$, $w$, and $w_2$. The energy at $H^2$ level becomes:
\be
\label{definehatE2}
\hat{\mathcal E_2} = \int |\eh_2|^2.
\ee
Moreover, with the bounds of the radiation \fref{radiation}, we can measure the difference between both energies at $H^2$ level.
\be
\label{chlocE2}
\hat{\mathcal E_2} \lesssim \mathcal E_2 + \int |H\zeta_b|^2 \lesssim \mathcal E_2 + b^2|\log b|^2.
\ee
Finally, $\hat w_2$ verifies the following partial differential equation:
\be
\pa_t \hat {w}_2 + \hat H_{\lambda} \hat{w}_2  =- \pa_t \tilde V \hat w + \hat H_{\lambda}\left( \frac{1}{\lambda^2} \hat{\mathcal F}_{\lambda} \right).
\ee
   \subsection{Modulation equations}
   With the choice of orthogonality conditions \fref{ortho}, we can now measure the error made taking $b = \frac{-\lambda_s}{\lambda}$ and $b_s = -b^2\left(1 + \frac 2{|\log b|}\right)$. This estimations are in the core of our proof. This demonstration is the same of in \cite{RS}, with the exception of a very small difference with the linear and non linear operators $L(\e)$ and $N(\e)$, which  truly brings any difficulty, because of the same interpolation bounds. In view of the importance to this lemma in our proof of the theorem \ref{main}, we considered useful to give integrally again this demonstration.
     \begin{lemma}[Modulation equations]
\label{modulationequations}
There holds the bound on the modulation parameters:
\be
\label{parameters}
\left|\frac{\lambda_s}{\lambda} + b\right| \lesssim \frac{b^2}{|\log b|} + \frac{1}{\sqrt{\log M}}\sqrt{\mathcal E_4},
\ee
\be
\label{parameterspresicely}
\left| b_s + b^2\left(1 + \frac{2}{|\log b|}\right) \right| \lesssim \frac{1}{\sqrt{\log M}} \left( \sqrt{ \mathcal E_4} +  \frac{b^2}{|\log b|}  \right).
\ee
\end{lemma}

\begin{remark} Note that this implies in the bootstrap the rough bounds:
\be
\label{rougboundpope}
|b_s|+\left|\lsl+b\right|\leq 2b^2.
\ee
and in particular \fref{bornes} holds. 
\end{remark}
{\bf Proof of Lemma \ref{modulationequations}}\\
\par
{\bf Step 1} Law for b.\\

Let
$$V(t) = |b_s + b^2| + \left|\frac{\lambda_s}{\lambda} + b\right|.$$
We take the inner product of \fref{eqofepsilon} with $H\Phi_M$ and estimate each terms.
\bee
&&\left(\pa_s \e,H \Phi_M\right) + \left(H \e,H \Phi_M\right) - \left( \frac{\lambda_s}{\lambda} \Lambda \e,H \Phi_M\right) \\ 
&=&-\left(\tilde{\Psi}_b,H \Phi_M\right) - \left(Mod(t),H \Phi_M\right) + \left(L(\e),H \Phi_M\right) + \left(N(\e),H \Phi_M\right).
\eee
First according to our choice of orthogonality \fref{ortho} 
\be
\left(\pa_s \e,H \Phi_M\right) = \pa_s(H\e,\Phi_M) = 0.
\ee
Then
\be
\left(H \e,H \Phi_M\right) = \left(H^2 \e,\Phi_M\right) \lesssim \|H^2 \e \|_{L^2} \| \Phi_M\|_{L^2} \lesssim \sqrt{\mathcal E_4 log M}.
\ee
From the construction of the profile, \fref{defmode} and the localization $\mbox{Supp}(\Phi_M)\subset[0,2M]$ from \fref{defdirection}:
 \bee
\nonumber \left(H\left(Mod(t)\right),\Phi_M\right) &=& -\left(b+\frac{\lambda_s}{\lambda}\right) \left(H\Lambda \qbt , \Phi_M\right)+ \left(b_s + b^2\right)\left(H\left(\tilde{T}_1 + 2b \tilde{T}_2\right),\Phi_M\right)\\
& = &  -(\Lambda Q, \Phi_M)(b_s + b^2)+O\left(c(M)b| V(t)|\right).
\eee
Using the Hardy bounds of Appendix B:
$$\left|\left(-\lsl\Lambda \e+L(\e)+N(\e),H\Phi_M\right)\right|\lesssim C(M)b(\mathcal E_4+|V(t)|).$$
Conclude from \fref{prodscalar} and the fundamental flux computation \fref{fluxcomputationonebis}:
\bee
b_s+b^2& = & \frac{(\Psit_b,H\Phi_M)}{(\Lambda Q,\Phi_M)}+O\left(\frac{\sqrt{\log M\mathcal E_4}}{\log M}\right)+O\left(C(M)b|V(t)|\right) \\
&=& -\frac{2b^2}{|\log b|}\left( 1 + O\left(\dfrac{1}{|\log b|}\right) \right) + O\left( \sqrt{\dfrac{\mathcal E_4}{\log M}}+C(M)b|V(t)|\right)
\eee
and \fref{parameterspresicely} is proved.
\\
\par
{\bf Step 2} Degeneracy of the law for $\lambda$.\\

 Now we take the inner product of \fref{eqofepsilon} with $\Phi_M$ and obtain:
$$(\mbox{Mod}(t),\Phi_M)=-(\Psit_b,\Phi_M)-(\pa_s \e + H \e,\Phi_M)-\left(-\lsl\Lambda \e+L(\e)+N(\e),\Phi_M\right).$$
From our choice of orthogonality conditions \fref{ortho}: $$(\pa_s \e + H\e,\Phi_M)=0.$$ 
In the same way as the last step, using the Hardy bounds of Appendix B:
$$\left|\left(-\lsl\Lambda \e+L(\e)+N(\e),\Phi_M\right)\right|\lesssim C(M)b\sqrt{\mathcal E_4}.$$
Next, we compute from  \fref{prodscalar} and the orthogonality \fref{intM}:
\bee
\nonumber \left(\mbox{Mod}(t),\Phi_M\right) &=& -\left(b+\frac{\lambda_s}{\lambda}\right) \left(\Lambda \qbt , \Phi_M\right)+ \left(b_s + b^2\right)\left(\tilde{T}_1 + 2b \tilde{T}_2,\Phi_M\right)\\
& = &  - 4 \log M (1+o_{M \to + \infty}(1))\left(\lsl+b\right)+O\left(C(M)b| V(t)|\right).
\eee
and observe the cancellation from \fref{sigmab0}, \fref{intM}:
\bee
\nonumber \left| \left(\Psit_b,\Phi_M \right) \right| \lesssim b^2|(\Sigma_b, \Phi_M)| + O(C(M)b^3)= c_b b^2 |(T_1,\Phi_M)| +O(C(M)b^3)=O(C(M)b^3).
\eee
We thus obtain the modulation equation for scaling:
\be
\label{pepepe}
\left|\lsl+b\right|\lesssim b^3C(M)+bC(M)O\left(\sqrt{\mathcal E_4}+ |V(t)|\right).
\ee
With \fref{parameterspresicely}, we obtain the bound $$|V(t)|\lesssim \frac{b^2}{|\log b|} + \frac{1}{\sqrt{log M}}\sqrt{\mathcal E_4}.$$
Injecting this bound in \fref{pepepe} implies the refined bound \fref{parameterspresicely}. This concludes the proof of Lemma \ref{modulationequations}.

    \subsection{Proof of the Theorem \ref{main}}
In this section, we conclude the proof of Theorem \ref{main} assuming the following a priori bounds on the solution on its maximum time interval of existence $[0,T)$, $0<T\leq +\infty$:
  \begin{itemize}

  \item{Energy estimates}
  \be
  \label{firstgg}
  \forall t \in [0,T[   \ \ \ \ \mathcal E_1(t) \leq \delta (b^*) \mbox{,}  \ \ \ \ \ \mathcal E_2(t) \lesssim b(t)^2 |\log b(t)|^5 \mbox{,}  \ \ \ \ \ \mathcal E_4(t) \lesssim \frac{b(t)^4}{|\log b(t)|^2}
  \ee
  \item{Link between both laws $b(t)$ and $\lambda(t)$}: there exist $ \alpha_1$, $\alpha_2 >0$ such that:
    \be
    \label{secondgg}
  C(u_0) b(t) |\log b(t)|^{\alpha_1} \leq \lambda(t) \leq C'(u_0) b(t) |\log b(t)|^{\alpha_2}
  \ee
  \end{itemize}
   
  \medskip
  \medskip  
  \par
  
  The heart of our analysis in section \ref{sectionfour} will be to produce such kind of solutions.\\
   We now assume \fref{firstgg} and \fref{secondgg} and prove Theorem \ref{main}. The proof adapts the argument in \cite{RS} which we sketch for the convenience of the reader.
   \\
   \par
{\bf Step 1} Finite time blow up.\\

Let $T\leq +\infty$ be the life time of $u$.  From \fref{firstgg}, \fref{secondgg}, $$-\frac{d}{dt}\sqrt\lambda=-\frac{1}{2\l\sqrt\l}\lsl\gtrsim \frac{b}{\l \sqrt{\lambda}}\gtrsim C(u_0)>0$$ and thus $\lambda$ touches zero in finite time which implies $$T<+\infty.$$ The bounds \fref{firstgg} and standard $H^4$ local well posedness theory ensure that blow up corresponds to
\be
\lambda(t) \rightarrow 0 \ \ \ \mbox{as} \ \ \ t \rightarrow T
\ee 
and thus, with \fref{secondgg}
\be
\label{boudnaryb}
\lambda(T)=b(T)=0.
\ee
\\
\par
{\bf Step 2} Derivation of the sharp blow up speed.\\

To begin, the modulation laws become with the bound \fref{firstgg} :
\be
\label{teteter}
\left| b + \frac{\lambda_s}{\lambda}\right| \lesssim \frac{b^2}{|\log b|}
\ee
\be
\label{parameterspresicelybis}
\left| b_s + b^2\left(1 + \frac{2}{|\log b|}\right) \right| \lesssim \frac{1}{\sqrt{\log M}} \frac{b^2}{|\log b|}.
\ee
We have defined M as a enough large constant. Let 
\be
\label{defbdelta}
B_{\delta} = \frac{1}{b^{\delta}}.
\ee
Since the variation of $b(s)$ is very small, we can consider in the first time that it's possible to take $M=B_{\delta}$ in \fref{parameterspresicelybis}. Assuming this, prove \fref{law}, and demonstrate afterwards that the made error is negligible for $\delta$ enough small. We have thus that \fref{parameterspresicelybis} becomes :
\be
\label{parameterspresicelyter}
\left| b_s + b^2\left(1 + \frac{2}{|\log b|}\right) \right| \lesssim \frac{1}{\sqrt{\log B_{\delta}}} \frac{b^2}{|\log b|} \lesssim \frac{b^2}{|\log b|^{\frac32}}.
\ee
We now integrate this in time using $\lim_{s\to +\infty}b(s)=0$ :
\be
\label{esitmateforb}
b(s)=\frac1s-\frac{2}{s\log s}+O\left(\frac{1}{s|\log s|^{\frac32}}\right).
\ee 
Using this decomposition of $b(s)$ in the modulation equation \fref{teteter}, we conclude: $$-\lsl=\frac1s-\frac{2}{s\log s}+O\left(\frac{1}{s|\log s|^{\frac32}}\right).$$ 
We rewrite this as $$\left|\frac{d}{ds}\log\left(\frac{s\lambda(s)}{(\log s)^2}\right)\right|\lesssim \frac{1}{s|\log s|^{\frac32}}$$ and thus integrating in time yields the existence of $\kappa(u_0)>0$ such that: 
$$\frac{s\lambda(s)}{(\log s)^2}=\frac{1}{\kappa(u_0)}\left[1+O\left(\frac{1}{|\log s|^{\frac32}}\right)\right].
$$
Taking the log yields the bound $$|\log \lambda|=|\log s|\left[1+O\left(\frac{|\log \log s|}{\log s}\right)\right]$$ and thus $$\frac{1}{s}=\kappa(u_0)\frac{\lambda}{|\log \lambda|^2}\left(1+o(1)\right).$$ Injecting this into \fref{esitmateforb} yields:
\be
\label{cnoencoenoe}
-\lambda\lambda_t=-\lsl=\frac1s\left(1+o(1)\right)=\kappa(u_0)\frac{\lambda}{|\log \lambda|^2}\left(1+o(1)\right)
\ee and thus $$-|\log \lambda|^2\lambda_t=\kappa(u_0)(1+o(1)).$$ Integrating from $t$ to $T$ with $\lambda(T)=0$ yields 
$$\lambda(t)=\kappa(u_0)\frac{T-t}{|\log (T-t)|^2}\left[1+o(1)\right],$$ and \fref{law} is proved. 
\\
\par
Prove now that the made error is negligible. Indeed, we take the inner product of \fref{eqofepsilon} with $H \chi_{B_\delta}\Lambda Q$ and obtain:
\bea
\label{vnoooejeojier}
\nonumber &&\frac{d}{ds}\left\{(H\e,\chi_{B_\delta}\Lambda Q)\right\}-(H\e,\pa_s\chi_{B_\delta}\Lambda Q)+\lsl(\chi_{B_\delta}\Lambda Q,H\Lambda\e)+(H^2\e,\chi_{B_\delta}\Lambda Q)\\
& = & \left(H\left[-\Psit_b+L(\e)-N(\e)-Mod\right],\chi_{B_\delta}\Lambda Q\right).
\eea
We must estimate all terms in this identity. First, for $\delta$ small enough, we have the rough bound:
$$|(H\e,\pa_s\chi_{B_\delta}\Lambda Q)|+|\lsl(\chi_{B_\delta}\Lambda Q,H\Lambda\e)|+|(H[L(\e)-N(\e)],\chi_{B_\delta}\Lambda Q)|\lesssim \frac{b}{b^{C\delta}}\sqrt{\mathcal E_4}\lesssim \frac{b^2}{|\log b|^2}.$$ 
For the linear term, we have immediately: $$|(H^2\e,\chi_{B_\delta}\Lambda Q)|\lesssim \sqrt{\mathcal E_4}\sqrt{|\log b|}\lesssim \frac{b^2}{\sqrt{|\log b|}}.$$
 The $\Psit_b$ term is computed from \fref{newpsib}:
$$(-H\Psit_b,\chi_{B_{\delta}}\Psi_b)=-b^2(H\Sigma_b,\chi_{B_{\delta}}\Psi_b)+O\left(\frac{b^3}{b^{C\delta}}\right)=b^2c_b(\Lambda Q,\chi_{B_{\delta}}\Lambda Q)+O\left(\frac{b^2}{|\log b|^2}\right).$$ 
From \fref{defmode}, we have the following estimate for the modulation term:
\bee
(-HMod,\chi_{B_\delta}\Lambda Q) & = & \left(\lsl+b\right)(H\Lambda \qbt,\chi_{B_\delta}\Lambda Q)-(b_s+b^2)(H(\tt_1+2b\tt_2),\chi_{B_\delta}\Lambda Q)\\
& = & (b_s+b^2)(\Lambda Q,\chi_{B_\delta}\Lambda Q)+O\left(\frac{b}{b^{C\delta}}\frac{b^2}{|\log b|}\right).
\eee
We now inject the estimates into \fref{vnoooejeojier} and obtain:
$$(b_s+b^2)(\Lambda Q,\chi_{B_\delta}\Lambda Q)=\frac{d}{ds}\left\{(H\e,\chi_{B_\delta}\Lambda Q)\right\}-c_bb^2(\Lambda Q,\chi_{B_{\delta}}\Lambda Q)+O\left(\frac{b^2}{\sqrt{|\log b|}}\right)$$ which we rewrite using \fref{cb} and an integration by parts in time:
\bea
\label{intito}
&&\frac{d}{ds}\left\{b-\frac{(H\e,\chi_{B_\delta}\Lambda Q)}{(\Lambda Q,\chi_{B_\delta}\Lambda Q)}\right\}+b^2\left(1+\frac{2}{|\log b|}\right)\\
\nonumber & = & O\left(\frac{b^2}{|\log b|^{\frac32}}\right)+(H\e,\chi_{B_\delta}\Lambda Q)\frac{(\Lambda Q,\pa_s\chi_{B_\delta}\Lambda Q)}{(\Lambda Q,\chi_{B_\delta}\Lambda Q)^2}.
\eea
We now estimate: 
$$
\left|(H\e,\chi_{B_\delta}\Lambda Q)\frac{(\Lambda Q,\pa_s\chi_{B_\delta}\Lambda Q)}{(\Lambda Q,\chi_{B_\delta}\Lambda Q)^2}\right|\lesssim \frac{\sqrt{\mathcal E_4}}{b^{C\delta}}\frac{|b_s|}{b}\lesssim \frac{b^3}{b^{C\delta}},$$
$$\left|\frac{(H\e,\chi_{B_\delta}\Lambda Q)}{(\Lambda Q,\chi_{B_\delta}\Lambda Q)}\right|\lesssim \frac{\sqrt{\mathcal E_4}}{b^{C\delta}}\lesssim \frac{b^2}{b^{C\delta}}.
$$
We inject these bounds into \fref{intito} and conclude that the difference between $b$ and $\tilde b$ is given by
\be
\label{lienbbtilde}
\tilde{b}=b-\frac{(H\e,\chi_{B_\delta}\Lambda Q)}{(\Lambda Q,\chi_{B_\delta}\Lambda Q)}=b+O\left(\frac{b^2}{|\log b|^2}\right)
\ee
satisfies the pointwise differential control:$$\left|\tilde{b}_s+\bt^2\left(1+\frac{2}{|\log \tilde{b}|}\right)\right|\lesssim \frac{\tilde{b}^2}{|\log \tilde{b}|^{\frac32}},$$ which we rewrite : $$\frac{\bt_s}{\bt^2\left(1+\frac{2}{|\log \bt|}\right)}+1=O\left(\frac{1}{|\log \bt|^{\frac32}}\right).$$ 
 We now integrate this in time using $\lim_{s\to +\infty}\bt(s)=0$ from \fref{boudnaryb}, \fref{lienbbtilde} and get:
$$\bt(s)=\frac1s-\frac{2}{s\log s}+O\left(\frac{1}{s|\log s|^{\frac32}}\right)$$ and thus from \fref{lienbbtilde}: 
\be
\label{esitmateforb}
b(s)=\frac1s-\frac{2}{s\log s}+O\left(\frac{1}{s|\log s|^{\frac32}}\right).
\ee 
This conclude the proof.
\\
\par
{\bf Step 3} Quantization of the focused energy.\\

We now turn to the proof of \fref{convustarb}, \fref{regularityustar} and adapt the strategy in \cite{MR5}. We shall need the following bound, which is a direct consequence of our construction and \fref{firstgg}:
 \be
 \label{vounds}
 \forall t\in [0,T), \ \ \|\Delta \tilde{u}(t,x)\|_{L^2}\leq C(v_0).
 \ee
where 
\be
\label{defvtilde}
\tilde u (t,x) = u(t,x) - \frac{1}{\lambda(t)} Q\left(\frac x {\lambda(t)}\right)
\ee
The regularity of $v(t,x)$ outside the origin is a standard consequence of parabolic regularity. Hence there exists $u^*\in \dot{H}^1$ such that $$\forall R>0, \ \ \nabla u(t)\to \nabla u^*\ \ \mbox{in} \ \ L^2(|x|\geq R)\ \ \mbox{as}\ \ t\to T.$$ Moreover, $v$ is $\dot{H^1}$ bounded by decrease of energy, and thus recalling the decomposition \fref{defvtilde} and the uniform bound \fref{vounds}:
$$ \nabla \tilde{u}(t)\to \nabla u^*\ \ \mbox{in} \ \ L^2\ \ \mbox{and}\ \ \Delta u^*\in L^2$$ which concludes the proof of \fref{convustarb}, \fref{regularityustar}. This concludes the proof of Theorem \ref{main}.

\section{Description of the initial data and bootstrap}
\label{sectionfour}
 The proof of the Theorem \ref{main} consists now in the demonstration of the existence of initial data, close to Q, whose the timing will be in agreement with the assumed bounds \fref{firstgg} and \fref{secondgg}.
 We have already seen the condition of smallness of b to assure the uniqueness of the decomposition, through the implicit function theorem:
 \be
 \label{init1}
 0 < b(0) < b^*(M) \ll 1.
 \ee 
 To be large regarding the constrains \fref{firstgg}, we fix the initial generous bounds namely:
\be
\label{init2energy}
|\mathcal E_1(0)|  \leq b(0)^2,
\ee
and
\be
\label{init2}
|\mathcal E_2(0)| + |\mathcal E_4(0)| \leq b(0)^{10}.
\ee
Moreover, there is a crucial difference compared to \cite{RS}. The linear operator $H$ posseses a negative direction $\psi$, source of instability, which can be harmful to the blow up dynamics, if we don't control this. Therefore, we manage this as in \cite{HR}. We note 
\be
\kappa(t) = \left(\e(t),\psi\right),
\ee
and
\be
a^+ = \kappa(0) =  \left(\e(0),\psi\right).
\ee
We impose that:
\be
\label{init3}
\left|a^+\right|  \leq \frac{2b(0)^{\frac52}}{|\log b(0)|}.
\ee
\par
 The propagation of regularity by the parabolic heat flow ensures that these estimates hold on some time interval $[0,t_1)$ together with the regularity $(\lambda, b)\in \mathcal C^1([0,t_1),\Bbb R^*_+\times \Bbb R)$. Given a large enough universal constant $K>0$ -independent of $M$-, we assume on $[0,t_1)$: 
\begin{itemize}
\item Control of b(t):
\be
\label{init1h}
0 <b(t)<10b(0).
\ee
\item Control of the radiation:
\be
\label{init2h}
\int | \nabla \varepsilon (t)|^2  \leq 10\sqrt{b(0)},
\ee
\be
\label{init3h}
|\mathcal E_2(t)| \leq K b^2(t)|\log b(t)|^5,
\ee
\be
\label{init3hbis}
|\mathcal E_4(t)| \leq K \frac{b^4(t)}{|\log b(t)|^2}.
\ee
\item \begin{it} A priori \end{it} bound on the unstable mode
\be
\label{init4h}
|\kappa(t)| \leq 2 \frac{b^{\frac52}}{|\log b|}.
\ee
\end{itemize} 
We may describe the bootstrap regime as follow :

\begin{definition}[Exit time]
\label{exittime}
Given $a^+ \in \left[ -2\dfrac{b(0)^{\frac52}}{|\log b(0)|};2\dfrac{b(0)^{\frac52}}{|\log b(0)|}\right]$, we let $T(a^+)$ be the life time of the solution to \fref{main} with initial data \fref{init1}, \fref{init2energy}, \fref{init2} and \fref{init3}, and $T_1(a^+) >0$ be the supremum of $T \in (0,T(a^+))$, such that for all $t \in [0;T]$, the estimates \fref{init1h}, \fref{init2h}, \fref{init3h}, \fref{init3hbis} and \fref{init4h} hold.
\end{definition}
The existence of blow up solutions in the regime described by Theorem \ref{main} now follows from the following:
\begin{proposition}
\label{bootstrapregime}
There exists $a^+ \in \left] -2\dfrac{b(0)^{\frac52}}{|\log b(0)|};2\dfrac{b(0)^{\frac52}}{|\log b(0)|}\right[$ such that
\be
T_1(a^+) = T(a^+)
\ee 
and then corresponding solution of \fref{main} blows up in finite time in the regime described by Theorem \ref{main}
\end{proposition}
We shall use the same strategy as in \cite{RR}, \cite{MRR}, \cite{HR}, and \cite{RS}. We will process in three times:
\begin{itemize}
\item First, we shall derive of suitable Lyapounov functionals at  Sobolev respectively $\dot H^4$ and $\dot H^2$  levels. That is the most difficult part of the proof, particularly because of the estimates of non linear terms,
  for whose we must  make a sharp study. Moreover, we shall see that it was crucial that $|a^+| \lesssim \frac{b(0)^{\frac52}}{|\log b(0)|}$.
\item
Secondly, we will reintegrate this functionals, to obtain improved bounds for $\mathcal E_2$ and $\mathcal E_4$. The bounds \fref{init1h} is a direct consequence of the energy decrease. Thus, only the last bounds \fref{init4h} for the unstable direction can be the cause of an exit time less than the life time of the solution.
\item Finally, we will study the dynamics of the unstable mode, and see that we can choose a $a^+$ to obtain Proposition \ref{bootstrapregime}. To ensure the existence of this solution, it is important that $|a^+| \gtrsim g(b(0))$ with $ \frac{b(x)^{3}}{|\log b(x)|} = o\left(g(x)\right)$ as $x \rightarrow 0$, hence the choice for the bounds of $a^+$  in \fref{init3}.
\end{itemize}
\section{Lyapounov monotonicities} 

\subsection{At $\dot H^4$ level}


\begin{proposition}[Lyapounov monotonicity $ \dot H^4$]
\label{AEI2}
There holds:
\bea
\label{monoenoiencle}
\frac{d}{dt} \left\{\frac{1}{\lambda^6}\left[\mathcal E_4+O\left(\sqrt b\frac{b^4}{|\log b|^2}\right)\right]\right\} \leq C\frac b {\lambda^{8}}\left[ \frac{\mathcal E_4}{\sqrt{\log M}}+\frac{b^4}{|\log b|^2}+\frac{b^2}{|\log b|}\sqrt{\mathcal E_4} \right]
\eea
for some universal constant $C>0$ independent of $M$ and of the bootstrap constant $K$ in \fref{init1h}, \fref{init2h}, \fref{init3h}, \fref{init3hbis}, provided $b^*(M)$ in \fref{init1} has been chosen small enough.
\end{proposition}
{\bf Proof of Proposition \ref{AEI2}}\\
\par
We recall the partial differential equations satisfied by $w_2$ and $w_4$:
\be
\pa_t w_2 + H_{\lambda} w_2 = - \pa_t \tilde{V} w +  H_{\lambda}\left(\frac1{\lambda^2}\mathcal{F}_{\lambda}\right) .
\ee
\be
\pa_t w_4 + H_{\lambda} w_4 = - \pa_t \tilde{V} w_2 - H_{\lambda}\left(\pa_t \tilde{V}w\right)  +  H^2_{\lambda}\left(\frac1{\lambda^2}\mathcal{F}_{\lambda}\right) .
\ee
Morever, we recall the action of time derivates on rescaling:
\bee
\label{actionderivatetime}
\pa_t v_{\lambda} (r) = \frac{1}{\lambda^2} \left( \pa_s v - \frac{\lambda_s}{\lambda}\Lambda v \right)_{\lambda}.
\eee
\par
{\bf Step 1} Energy identity
\begin{lemma} [Energy identity   $\dot H^4$]
\label{energyH4}
\bee
&&\frac12 \frac d {dt} \int \left\{w_4^2 -2 \pa_{t} \tilde{V} ww_4\right\} \\
&=& - \int w_4 H_{\lambda}w_4 + \int \left( \pa_t \tilde{V}\right)^2 w_2w - \int  \pa_{tt} \tilde{V} ww_4 +  \int H_{\lambda}\left( \pa_t \tilde{V} w \right)  \pa_t \tilde{V} w \\
&+& \int w_4 H^2_{\lambda}\frac1{\lambda^2}\mathcal{F}_{\lambda}  - \int H_{\lambda}\left( \pa_t \tilde{V} w \right) H_{\lambda}\left(\frac1{\lambda^2}\mathcal{F}_{\lambda}\right) -\int \pa_t \tilde{V} w_4\frac1{\lambda^2}\mathcal{F}_{\lambda}.
\eee
\end{lemma}
{\it Proof of the Lemma \ref{energyH4}:} We propose here a simplification with respect to the algebra in \cite{RS}. Dissipation also allows us to sign some terms and avoid the study of suitable quadratic forms as in \cite{HR}. We compute the energy identity:
\\
\\
\bee
\frac12 \frac d {dt} \int w_4^2 &=& \int w_4 \pa_t w_4\\
&=& \int w_4 \left( -H_{\lambda}w_4 - \pa_t \tilde{V} w_2 - H_{\lambda}\left( \pa_t \tilde{V} w \right)+  H^2_{\lambda}\frac1{\lambda^2}\mathcal{F}_{\lambda} \right) .
\eee  
We now treat separately the second and the third term. 
\bee
  &&- \int w_4 \pa_t \tilde{V} w_2 \\
  &=&  \int \pa_t \tilde{V} w_2 \left( \pa_t w_2  + \pa_t \tilde{V} w -  H_{\lambda}\left(\frac1{\lambda^2}\mathcal{F}_{\lambda}\right) \right)\\
&=&  \int \left( \pa_t \tilde{V}\right)^2 w_2w - \int \pa_t \tilde{V} w_2 H_{\lambda}\left(\frac1{\lambda^2}\mathcal{F}_{\lambda}\right) + \int  \pa_t \tilde{V}  \pa_t \left( \frac{w_2^2}{2}\right)\\
&=&  \int \left( \pa_t \tilde{V}\right)^2 w_2w - \int \pa_t \tilde{V} w_2 H_{\lambda}\left(\frac1{\lambda^2}\mathcal{F}_{\lambda}\right) - \frac{1}{2}\int  \pa_{tt} \tilde{V} w_2^2 + \frac{1}{2}\frac{d}{dt}\left(\int  \pa_{t} \tilde{V} w_2^2\right).
\eee
Now
\bee
-\int w_4  H_{\lambda}\left( \pa_t \tilde{V} w \right) &=&  \int H_{\lambda}\left( \pa_t \tilde{V} w \right) \left( \pa_t w_2 + \pa_t \tilde{V} w -  H_{\lambda}\left(\frac1{\lambda^2}\mathcal{F}_{\lambda}\right) \right)\\
&=&  \int H_{\lambda}\left( \pa_t \tilde{V} w \right) \left(\pa_t \tilde{V} w -  H_{\lambda}\left(\frac1{\lambda^2}\mathcal{F}_{\lambda}\right) \right)\\ &+& \frac{d}{dt} \int H_{\lambda}\left( \pa_t \tilde{V} w \right) w_2 -  \int \pa_t\left[ H_{\lambda}\left( \pa_t \tilde{V} w \right) \right]w_2.\\
\eee
The last term becomes
\bee
-  \int \pa_t H_{\lambda}\left( \pa_t \tilde{V} w \right) w_2 &=& \int \left( \pa_t \tilde{V}\right)^2 w_2w - \int \pa_{tt} \tilde{V} w w_4 - \int \pa_t \tilde{V} \pa_tw w_4,
\eee
and
\bee
 - \int \pa_t \tilde{V} \pa_tw w_4 &=& \int \pa_t \tilde{V} \left(w_2 - \frac1{\lambda^2}\mathcal{F}_{\lambda} \right)w_4\\
 &=&- \int \pa_t \tilde{V} \frac1{\lambda^2}\mathcal{F}_{\lambda}w_4 + \int \pa_t \tilde{V} w_2 \left(- \pa_t w_2 - \pa_t \tilde{V} w +  H_{\lambda}\left(\frac1{\lambda^2}\mathcal{F}_{\lambda}\right) \right) \\
 &=&- \int \pa_t \tilde{V} \frac1{\lambda^2}\mathcal{F}_{\lambda}w_4 + \int \pa_t \tilde{V} w_2 \left( - \pa_t \tilde{V} w +  H_{\lambda}\left(\frac1{\lambda^2}\mathcal{F}_{\lambda}\right) \right)  \\
 &-& \frac12\frac{d}{dt} \left(\int \pa_t \tilde{V}w_2^2\right) + \frac12 \int \pa_{tt}\tilde{V}w_2^2 .
\eee
In the following steps, we estimate each terms of Lemma \ref{energyH4} in order to prove Proposition \ref{AEI2}.\\
\par 
{\bf Step 2} Lower order quadratic terms
\\
\\
We have from \fref{comportementv} and \fref{actionderivatetime}, and the modulation equations, the bounds:
\be
\label{dtV}
|\pa_t \tilde V|\lesssim \frac{b}{\lambda^4}\frac{1}{1+y^4} \ \ \ \ |\pa_{tt} \tilde V|\lesssim \frac{b}{\lambda^6}\frac{1}{1+y^4}.
\ee

Using \fref{subpositivity}, we obtain
\bee
-\int H_{\lambda}w_4w_4 = - \frac1{\lambda^8}\int H \e_4\e_4 &\lesssim& \frac1{\lambda^8}\left( \int \e_4 \psi\right)^2\lesssim \frac{\zeta^2}{\lambda^8} \left(\int \e \psi \right)^2 \\
&\lesssim& \frac1{\lambda^8}\kappa^2 \lesssim \frac b{\lambda^8}\frac{b^4}{|\log b|^2}.
\eee
Remark that:
\bee
\int H_{\lambda}\left( \pa_t \tilde{V} w \right)  \pa_t \tilde{V} w =  \int \left( \pa_t \tilde{V}\right)^2 w_2w - \int \Delta\left( \pa_t \tilde V\right) \pa_t \tilde V w^2 - 2 \int \pa_{rt} \tilde V \pa_t \tilde V \pa_r w w.
\eee
We treat now the two following terms:
\bee
&&\left| 2 \int \left( \pa_t \tilde{V}\right)^2 w_2w - \int  \pa_{tt} \tilde{V} ww_4 - \int \Delta\left( \pa_t \tilde V\right) \pa_t \tilde V w^2 - 2 \int \pa_{rt} \tilde V \pa_t \tilde V \pa_r w w \right| \\
&\lesssim& \frac{b^2}{\lambda^8} \int \left( \frac{\e \e_2}{1+y^8} + \frac{\e\e_4}{1+y^4} + \frac{\e^2}{1+y^{10}} + \frac{\e \pa_y \e}{1+y^9}\right) \lesssim \frac{b^2}{\lambda^8} b^4 |\log b|^C.
\eee
The last inequality comes from Cauchy-Schwarz and the bounds \fref{estun} and \fref{lossyboundwperp}.
Finally, we estimate the boundary term in time
\bee
\left |\int \pa_{t} \tilde{V} ww_4\right| \lesssim \frac{b}{\lambda^6}\left(\int \frac{\e^2}{1+y^8}\right)^{\frac12}\left( \int \e_4^2\right)^{\frac12} \lesssim \frac{b}{\lambda^6}b^4 |\log b|^C\lesssim \frac{\sqrt b}{\lambda^6}\frac{b^4}{|\log b|^2}.
\eee
\\ \par
{\bf Step 3} Further use of dissipation
\\
\\
First, we claim the following bounds:
\bea
\label{y8bound} \int \frac{1}{1+y^8}\mathcal F^2 &\lesssim& \left[\frac{b^4}{|\log b|^2} + \frac{\mathcal E_4}{\log M}\right] \\
\label{H2boundF} \int |H^2 \mathcal F|^2 &\lesssim&  b^2 \left[\frac{b^4}{|\log b|^2}+ \frac{\mathcal E_4}{\log M}\right].
\eea
Thus,
\bee
&&\left| \int w_4 H^2_{\lambda}\frac1{\lambda^2}\mathcal{F}_{\lambda}  - \int H_{\lambda}\left( \pa_t \tilde{V} w \right) H_{\lambda}\left(\frac1{\lambda^2}\mathcal{F}_{\lambda}\right) -\int \pa_t \tilde{V} w_4\frac1{\lambda^2}\mathcal{F}_{\lambda}\right|\\
& \lesssim&\frac{1}{\lambda^8} \left\{\int \left(\e_4^2 + b^2\frac{\e^2}{1+y^8} \right)^{\frac12}\left( \int |H^2 \mathcal F|^2\right)^{\frac12} + b \left(\int\e_4^2 \right)^{\frac12}\left( \int \frac{1}{1+y^8}\mathcal F^2\right)^{\frac12} \right\} \\
&\lesssim&\frac b {\lambda^{8}}\left[ \frac{\mathcal E_4}{\sqrt{\log M}}+\frac{b^4}{|\log b|^2}+\frac{b^2}{|\log b|}\sqrt{\mathcal E_4} \right].
\eee
This concludes the proof of the Proposition \ref{AEI2}. We now turn of the proof of \fref{y8bound} and \fref{H2boundF}. We recall that
\bee
\mathcal F = - \tilde{\Psi}_b- Mod(t)  + L(\e) + N(\e) .
\eee
\par
{\bf Step 4} $\tilde \Psi_b$ terms.
\\
\\
The contribution of the $\tilde \Psi_b$ terms in \fref{y8bound} and \fref{H2boundF} has already been proved in the Lemma \ref{localisation}. For that matter, the construction of the approximated solution by the profiles $(T_i)_{1 \leq i \leq 3}$ has been made to obtain these good estimates.
\\
\par
{\bf Step 5} $Mod(t)$ terms.
\\
\\
Recall the definition \fref{defmode} of $Mod(t)$.
\bee
Mod(t)=- \left( \frac{\lambda_s}{\lambda} +b\right) \Lambda \qbt + (b_s + b^2)(\tilde{T_1} +2b\tilde{T_2}).
\eee
With the modulation equation \fref{parameters} and \fref{parameterspresicely}, we have:
\bee
\left|\frac{\lambda_s}{\lambda} +b \right|^2 + \left| b_s + b^2\right|^2 \lesssim \frac{b^4}{|\log b|^2}+ \frac{\mathcal E_4}{\log M}.
\eee
But
\bee
\int \frac{1}{1+y^8}|\Lambda \tilde Q_b|^2 + \int \frac{1}{1+y^8}|\tilde{T_1} +2b\tilde{T_2}|^2 \lesssim 1.
\eee
Now,
\be
\label{estimationdernierarticle}
\int  |H^2\Lambda \tilde Q_b|^2 + \int |H^2(\tilde{T_1} +2b\tilde{T_2})|^2 \lesssim b^2.
\ee
Indeed :
\be
\int \left| H^2\tilde{T_1}\right|^2 \lesssim \int_{B_1 \leq y \leq 2B_1} \left| \frac{\log y}{y^4}\right| \lesssim \frac{|\log b|^2}{B_1^4} \lesssim b^2
\ee
There is a whole proof of the estimate for $\tilde{T_2}$ in \cite{RS}. Here is a summary of this demonstration. With the definition of $\tilde{T}_2$, and the bound \fref{decayttwo} of $T_2$, we have :
\be
\label{cpeopeepi}
\int|H^2\tt_2|^2\lesssim \left[\int_{B_1\leq y\leq 2B_1}\left|\frac{y^2}{y^4}\right|^2+\int_{y\leq 2B_1}|H\Sigma_2|^2\right]\lesssim 1+\int_{y\leq 2B_1}|H\Sigma_2|^2.
\ee From the construction of the radiation $\Sigma_b$ and the definition of $\Sigma_2$, we compute:
\bee
H\Sigma_2 & = &  H\Sigma_b + H(T_1 - \Lambda T_1) +O\left(\frac{y|\log y|^2}{1+y^5}\right)\\
& = & \frac{1}{|\log b|}O\left(\frac{1}{1+y^2}{\bf1}_{y\leq 3B_0}\right)+H(T_1 - \Lambda T_1) +O\left(\frac{y|\log y|^2}{1+y^5}\right).
\eee
But
\bee
HT_1-H\Lambda T_1 & = & HT_1-\left(2HT_1+\Lambda HT_1-\frac{\Lambda V}{y^2}T_1\right) \\
& = & O\left(\frac{\log y}{1+y^3}\right).
\eee
We thus conclude: $$\int |H\Sigma_2|^2\lesssim \frac{1}{|\log b|^2}\int_{y\leq 2B_0}\frac{1}{1+y^4}+\int_{y\leq 2B_1}\frac{|\log y|^4}{1+y^8}\lesssim 1,$$
and the contribution of $Mod(t)$ terms to \fref{y8bound} and \fref{H2boundF} are small enough.
\\
\par
{\bf Step 6} Small linear term $L(\e)$.
\\
\\
We recall the expression of $L(\e)$:
\bee
L(\e) = 3\left( \tilde Q_b^2 - Q^2\right) \e.
\eee
We have, with the rough bounds \fref{asympT1}, \fref{roughT2} and \fref{roughT3}:
\bee
\left|3\left( \tilde Q_b^2 - Q^2\right)\right| \lesssim  b
\eee
and, with \fref{lossyboundwperp}
\bee
\int \frac{1}{1+y^8}|L(\e)|^2 \lesssim b^2\int \frac{1}{1+y^8}\e^2 \lesssim \frac{b^4}{|\log b|^2}.
\eee
Let us study the second estimate. In order to do that, let
\be
g = 3\left( \tilde Q_b^2 - Q^2\right).
\ee
We have the following bound:
\be
\label{boundg}
|H^j\left(\pa^i_y g \right)| \lesssim b \frac{y^k(1+|\log y|)}{1+y^{2+i+2j+k}}, \ 0 \leq i \leq 4, \ 0 \leq j \leq 1
\ee
where
\bee
k = \max \{0;2-i-2j\}.
\eee
So, with the bounds \fref{boundg} and those of the Lemma \ref{lemmainterpolation}, we obtain:
\bee
&&\int |H^2(L(\e))|^2 =\int |H^2(g \e)|^2 = \int |H\left( g H\e - \e\Delta g - \pa_y g \pa_y \e \right)|^2 \\
&\lesssim&  \int |g H^2(\e)|^2 + \int |\Delta g H\e|^2 + \int |\pa_y g \pa_y( H\e)|^2 \\
&+&  \int |\Delta g H \e|^2 + \int |\e \Delta^2g |^2 + \int |\pa_y(\Delta g) \pa_y \e|^2 \\
&+&  \int |H(\pa_y g) \pa_y \e|^2 + \int |\pa_y g \Delta(\pa_y \e)|^2 + \int |\pa_{yy} g \pa_{yy}\e|^2 \\
&\lesssim& \frac{b^6}{|\log b|^2}.
\eee
In order to improve the redaction, we didn't have developed the sharp estimation for each terms. The method to use is the same one that we 	are going to use now for the nonlinear terms. 
\\ \par
{\bf Step 7} Nonlinear term $N(\e)$.
\\
\\
We recall the expression of $N(\e)$:
\bee
N(\varepsilon) = 3\tilde Q_b \varepsilon^2 + \varepsilon ^3.
\eee
We have, with the rough bounds \fref{asympT1}, \fref{roughT2} and \fref{roughT3}:
\bee
\left|3\left( \tilde Q_b(y) - Q(y)\right)\right| \lesssim  |b|(1+ |\log y|)
\eee
and, with \fref{epsilonborne}
\bee
&&\int \frac{1}{1+y^8}\left|3(\tilde Q_b - Q + Q) \varepsilon^2\right|^2 \lesssim b^2\int \frac{1+|\log y|^2}{1+y^8}\e^4 + \int \frac{\e^4}{1+y^{12}}
\\&\lesssim& b^2\left \| \e \right \|^4_{L^{\infty}} \int \frac{1+|\log y|^2}{1+y^8} + \left \| \frac{\e}{1+y^2}\right \|_{L^{\infty}}^2\int \frac{\e^2}{1+y^8}
\lesssim \frac{b^4}{|\log b|^2}
\eee
and to conclude:
\bee
&&\int \frac{1}{1+y^8}\left|\varepsilon^3\right|^2 
\lesssim \left \| \e \right \|^6_{L^{\infty}}\int \frac{1}{1+y^8} \lesssim \frac{b^4}{|\log b|^2}.
\eee
For the second bound, let us compute $H(\e^3)$
\bee
H(\e^3)&=&\e^2H(\e) - 2 \pa_y(\e^2)\pa_y \e - \Delta(\e^2)\e\\
&=&\e^2H(\e) - 2 (\pa_y \e)^2 \e - \left( \e\pa_{yy}\e  +(\pa_y \e )^2 + 3\frac{\e\pa_y \e}{y} \right)\e \\
&=&\e^2H(\e) - 3 (\pa_y \e)^2 \e - \e^2\Delta \e \\
&=& 2 \e^2H(\e) + V \e^3 - 3 \e(\pa_y \e)^2.
\eee
Now, we treat each terms separately. First:
\bee
&&\int |H(\e^2H(\e))|^2 \lesssim \int \e^4|H^2(\e)|^2+ \int \Delta(\e^2)^2H(\e)^2 + \int |\pa_y(\e^2)|^2|\pa_yH(\e)|^2 \\
&\lesssim& \left \| \e \right \|_{L^{\infty}}^4 \mathcal E_4 + \left(\left \| \e \right \|_{L^{\infty}_{y \geq 1}}^2\left \| \pa_{yy}\e \right \|_{L^{\infty}_{y \geq 1}}^2 + \left \| \pa_{y}\e \right \|_{L^{\infty}}^4 + \left \| \pa_y\e \right \|_{L^{\infty}_{y \geq 1}}^2\left \| \frac{\e}{y} \right \|_{L^{\infty}_{y \geq 1}}^2 \right) \mathcal E_2 \\
 &+& \left(\left \| \e \right \|_{L^{\infty}_{y \leq 1}}^2\left \| y\pa_{yy}\e \right \|_{L^{\infty}_{y \leq 1}}^2  + \left \| \pa_y\e \right \|_{L^{\infty}_{y \leq 1}}^2\left \| \e \right \|_{L^{\infty}_{y \leq 1}}^2 +  \left \| \e (1+|\log y|^2) \right \|_{L^{\infty}_{y \geq 1}}^2\left \| y \pa_y \e \right \|_{L^{\infty}_{y \geq 1}}^2 \right) \mathcal E_4 \\
&\lesssim& \frac{b^6}{|\log b|^2}.
\eee
Secondly, using that 
\be
\label{boundV}
|\pa^i_y V| \lesssim \frac{1}{1+y^{4+i}} \  \  \  0 \leq i \leq 2,
\ee
we have:
\bee
&& \int |H(V \e^3)|^2 \lesssim \int |H(V)|^2\e^6 + \int \Delta (\e^3)^2V^2+ \int |\pa_y(\e^3)|^2|\pa_y V|^2 \\
&\lesssim&  \left\| \e \right \|_{L^{\infty}}^2\left(\left\| \frac{\e}{1+y} \right \|_{L^{\infty}}^2\left \| \frac{\e}{1+y^2} \right \|_{L^{\infty}}^2 + \left \| \pa_{y}\e \right \|_{L^{\infty}}^4+ \left \| \frac{\e}{1+y^2} \right \|_{L^{\infty}}^2\left \| \pa_y\e \right \|_{L^{\infty}}^2  \right)\int \frac{1}{1+y^6} \\
&+&  \left\| \e \right \|_{L^{\infty}}^4 \int \frac{|\pa_{yy} \e|^2}{1+y^8} \lesssim \frac{b^6}{|\log b|^2}.
\eee
Lastly
\bee
&& \int |H(\e(\pa_y \e)^2)|^2 \lesssim \int  |H(\e)|^2|(\pa_y \e)^2)|^2+ \int \e^2\left| \Delta((\pa_y \e)^2)\right|^2+ \int \pa_y\e^2 \left|\pa_y((\pa_y \e)^2)\right|^2\\
&\lesssim& \left \| \pa_{y}\e \right \|_{L^{\infty}}^4 \mathcal E_2 + \int \e^2\left| \Delta((\pa_y \e)^2)\right|^2 \lesssim \frac{b^6}{|\log b|^2} + \int \e^2\left| \Delta((\pa_y \e)^2)\right|^2.
\eee
But
\bee
 &&\int \e^2\left| \Delta((\pa_y \e)^2)\right|^2 \\
 &=& \int \e^2 \left( |\pa_{yy} \e|^2 + \pa_y\e \pa_y^3 \e + 3 \frac{\pa_y \e \pa_{yy} \e}{y}\right)^2\\
 &=& \int \e^2 \left( \pa_{yy} \e \left(H(\e) + V \e\right)  + \pa_y\e \pa_y^3 \e \right)^2\\
 &\lesssim&  \left \| \e \right \|_{L^{\infty}}^2 \left( \left \| y\pa_{yy} \e \right \|_{L^{\infty}_{y \leq 1}}^2\mathcal E_4 +\left \| \pa_{yy} \e \right \|_{L^{\infty}_{y \geq 1}}^2\mathcal E_2 +  \left \| \e \right \|_{L^{\infty}}^2 \int \frac{|\pa_{yy}\e|^2}{1+y^8}\right) \\
 &+& \left \| \e (1+|\log y|^2) \right \|_{L^{\infty}_{y \geq 1}}^2\left \| y \pa_y \e \right \|_{L^{\infty}_{y \geq 1}}^2 \mathcal E_4  +\left \| \e \right \|_{L^{\infty}_{y \leq 1}}^2\left \|\pa_y \e \right \|_{L^{\infty}_{y \leq 1}}^2 \int _{y \leq 1}  |\pa^3_y \e|^2\\
 &\lesssim& \frac{b^6}{|\log b|^2} .
\eee
Thus,
\be
\label{firstpart}
\int |H^2(\e^3)|^2 \lesssim \frac{b^6}{|\log y|^2}.
\ee
Treat now the other contribution of $N(\e)$ in the bound \fref{H2boundF}. Let
\bee
f = 3 \left( \tilde Q_b - Q\right).
\eee
We have the following bounds:
\bea
\label{boundf}|\pa_y^i f | &\lesssim& b \frac{y^{2-i}(1+|\log y|)}{1+y^2} + \frac{1}{1+y^{2+i}}, \ \ \ 0 \leq i \leq 2 \\
\label{boundHf}|H(\pa_y^j f)| &\lesssim& b \frac{y^k}{1+y^{2+j+k}} + \frac{1}{1+y^{4+j}}, \ \ \ 0 \leq j \leq 2
\eea
where 
\bee
k= \left\{\begin{array}{ll}1\ \ \mbox{for}\ \ j= 1,\\ 0\ \ \mbox{otherwise}.\end{array}\right.
\eee
Let us compute $H(f\e^2)$
\bea
\nonumber H(f\e^2) &=& H(\e^2) f - \e^2\Delta f - 2 \pa_y f \pa_y (\e^2)\\
\label{decompositionHFE2} &=& H(\e^2) f +\e^2Hf - 2 \pa_y f \pa_y (\e^2) +Vf\e^2
\eea
In the same way as the last proof, we treat each term separately. First:
\be
\label{firstterm}
H\left(H(\e^2) f \right) = H^2(\e^2)f - \Delta f H(\e^2) - \pa_y f \pa_y H(\e^2).
\ee
Let us estimate the three components using that:
\be
\label{HE2}
H(\e^2) = 3 \e H(\e) + 2 V \e^2 - 2 \left( \pa_y \e\right)^2
\ee
and,
\bea
\label{H2E2}
H^2(\e^2) &=& 3\left(\e H^2 (\e) + |H(\e)|^2 + V\e H(\e) - 2 \pa_y \e \pa_y H(\e)\right) \\
\nonumber&+& 2\left(HV \e^2 - V \Delta(\e^2) - 2 \pa_y V \pa_y (\e^2) \right)\\ 
\nonumber&-& 2\left( 3\pa_y \e H(\pa_y \e) + 2 V (\pa_y \e)^2 - 2 \left( \pa_{yy} \e \right)^2 \right)
\eea
 and moreover,
\be
\label{HPAYE}
H(\pa_y \e) = \pa_yH(\e) - \frac{\pa_y \e}{y^2} - \pa_y V \e
\ee
\fref{firstterm}, \fref{HE2}, \fref{H2E2} and \fref{HPAYE} together with the bounds \fref{boundf}, \fref{boundHf}, \fref{boundV} and those of the Lemma \ref{lemmainterpolation} imply:
\be
\int |H\left(H(\e^2) f \right) |^2 \lesssim \frac{b^6}{|\log b|^2}
\ee
Let us study the second term of \fref{decompositionHFE2} .
\be
\label{HE2HF}
H\left(\e^2Hf \right) = \e^2H^2(f) - Hf \Delta(\e^2) - 2 \pa_y (Hf) \pa_y (\e^2).
\ee
The bounds \fref{boundf}, \fref{boundHf}, the Lemma \ref{lemmainterpolation} yields 
\be
\int |H\left(\e^2Hf \right) |^2 \lesssim \frac{b^6}{|\log b|^2}.
\ee
We estimate the two last terms in \fref{decompositionHFE2} in the same way. This concludes the proof  of \fref{H2boundF} and thus of Proposition \ref{AEI2}.
\subsection{At $\dot H^2$ level}
For the $H^2$ level, we use the profile $\hat Q_b $ localized near $B_0$. The description of this localization and the estimates of the new error generated by these are given by the Lemma \ref{localisationbis} and in the subsection \ref{eed}. We recall the equation verified by $\hat{w}_2$:
\be
\pa_t \hat {w}_2 + \hat H_{\lambda} \hat{w}_2  =- \pa_t \tilde V \hat w + \hat H_{\lambda}\left( \frac{1}{\lambda^2} \hat{\mathcal F}_{\lambda} \right).
\ee
\par 
\begin{proposition}[Lyapounov monotonicity $\dot H^2$]
\label{ttpptt}
There holds:
\be
\label{monotonicityE2}
\frac{d}{dt} \frac{\hat{\mathcal E}_2}{\lambda^2}\lesssim \frac{b^3|\log b|^2}{\lambda^4}
\ee
\end{proposition}
{\bf Proof of the Proposition \ref{ttpptt}}
\\
\bea
\nonumber \frac12 \frac d {dt} \frac{\hat{\mathcal E}_2}{\lambda^2} &=& \int \hat{w}_2\left[  - H_{\lambda} \hat{w}_2  - \pa_t \tilde V \hat w + H_{\lambda}\left( \frac{1}{\lambda^2} \hat{\mathcal F}_{\lambda} \right)\right]\\
&=& - \int \hat{w}_2H_{\lambda}\hat{w}_2 - \int \pa_t \tilde V \hat w \hat w_2 + \int \hat w_2 H_{\lambda}\left( \frac{1}{\lambda^2} \hat{\mathcal F}_{\lambda} \right). 
\eea
We study each terms separately:
\\
\par
To begin, in agreement with \fref{subpositivity}, the \begin{it}a priori\end{it} bound \fref{init4h} on the unstable direction, and the bounds \fref{radiationbis} we have
\bea
\nonumber - \int \hat{w}_2H_{\lambda}\hat{w}_2 &\lesssim& \frac{1}{\lambda^4}(\hat {\e}_2,\psi)^2 \\
\nonumber &\lesssim& \frac{1}{\lambda^4}\left[ (\e_2,\psi)^2 + (H\zeta_b,\psi)^2\right] \\
\label{aaa1}&\lesssim& \frac{1}{\lambda^4}\left[ \kappa^2 + \|H^2\zeta_b\|_{L^2}^2\right] \lesssim \frac{b^3|\log b|^2}{\lambda^4}.
\eea
\\
The second term is a lower order quadratic term. Recall that:
\bee
|\pa_t \tilde V| \lesssim \frac{b}{\lambda^4} \frac{1}{1+y^4}.
\eee
Hence, using Cauchy-Schwarz, with the \begin{it}a priori\end{it} bound \fref{init3h}, \fref{chlocE2} measuring the difference between the two energies at $\dot H^2$ level, and the bounds \fref{radiationbis} and \fref{estun} 
\be
\label{aaa2}
\left|\int \pa_t \tilde V \hat w \hat w_2\right| \lesssim  \frac{b}{\lambda^4} \hat{\mathcal E}_2^\frac12\left(\int \frac{ |\e|^2 + |\zeta_b|^2}{1+y^8} \right)^\frac12 \lesssim \frac{b}{\lambda^4} b^3 |\log b|^C.
\ee
\\
For the last, we don't use exactely the same strategy as for the control $\dot H^4$. Indeed, for the term of error, the term of modulation, the global $L^2$ bounds for $\hat {\e}_2$, that we dispose is too rough. We must then improve it, as both terms are localized for $y \leq 2B_0$. So:
\bea
\label{improvedboundE2}
\nonumber \int_{y\leq 2B_0}|\eh_2|^2 & \lesssim & B_0^4|\log b|^2\int \frac{|\e|^2}{(1+y^4)|\log y|^2}+\int|H\zeta_b|^2\\
& \lesssim & C(M)b^2+b^2|\log b|^2\lesssim b^2|\log b|^2.
\eea
The term of error $\Psih_b$ is now estimated using  \fref{controleh2erreurtildehat} and the improved bound \fref{improvedboundE2}:
\be
\label{aaa3}
|(\eh_2,H\Psih_b)| \lesssim  \|H\Psih_b\|_{L^2}\|\eh_2\|_{L^2(y\leq 2B_0)}\lesssim  \left(b^4|\log b|^2b^2|\log b|^2\right)^{\frac12}\lesssim b^3|\log b|^2.
\ee
 We next estimate from \fref{asympT1}, \fref{decayttwo}:
$$ \int|H\hat{T}_1|^2\lesssim \int_{y\leq 2B_0}|\Lambda Q|^2+\int_{B_0\leq y\leq 2B_0}\left|\frac{\log y}{y}\right|^2\lesssim |\log b|^2,$$
$$\int |H\hat{T}_2|^2\lesssim  \int_{y\leq 2B_0}\left|\frac{y}{y^2b|\log b|}\right|^2\lesssim \frac{1}{b^2|\log b|},$$ and thus from \fref{defmodebis}, \fref{parameters}, \fref{parameterspresicely}: 
\bee
\int|H\widehat{Mod}(t)|^2 & \lesssim & \left|\lsl+b\right|^2\int|H\Lambda \qbh|^2+|b_s+b^2|^2\int|H(\hat{T}_1+b\hat{T}_2)|^2\\
& \lesssim & \frac{b^4}{|\log b|^2}|\log b|^2\lesssim b^4|\log b|^2.
\eee
Moreover, $\mbox{Supp}(H\widehat{Mod})\subset[0,2B_0]$ and thus with \fref{improvedboundE2}:
\be
\label{aaa4}
|(\eh_2,H\widehat{Mod})| \lesssim  \left(b^4|\log b|^2b^2|\log b|^2\right)^{\frac12}\lesssim b^3|\log b|^2.
\ee
We now claim the following bound for the small linear term, and the non linear term:
\be
\label{HhatLHhatN}
\int |H\hat L (\hat {\e})|^2 + |H\hat N (\hat {\e})|^2 \lesssim b^5.
\ee
Assume \fref{HhatLHhatN}. Thus,
\be
\label{aaa5}
|(\eh_2,HL (\hat {\e}))| \lesssim b |\log b|^C b^{\frac52} \lesssim b^3 |\log b|^2.
\ee
\fref{aaa5}, together with \fref{aaa1},  \fref{aaa2},  \fref{aaa3} and  \fref{aaa4} concludes the proof of the Proposition \ref{ttpptt}. 
\\
\\
{\it Proof of \fref{HhatLHhatN}}: 
\\
\\
We recall that:
\bee
\hat L(\eh) = 3(\hat Q_b^2 - Q^2) \eh.
\eee
In the same way as previously, we let:
\bee
\hat g = 3(\hat Q_b^2 - Q^2),
\eee
for which we have the bounds:
\bee
| \hat g| + |\pa_y \hat g| + |H \hat g  | &\lesssim& b \ \ \ \mbox{for} \ \ \ y \leq 1 \\
\left| \frac{\hat g}{1+|\log y|} \right|+ \left|y\pa_y \hat g \right|+ \left|y^2H \hat g  \right| &\lesssim& \frac b{1+y^2} \ \ \ \mbox{for} \ \ \ y \geq 1.
\eee 
This bounds, and \fref{radiation}, \fref{radiationbis}, and whose of the Lemma \ref{lemmainterpolation} together with the following decomposition:
\bee
H \hat L = \eh H \hat g - \hat g \Delta \eh - 2 \pa_y \hat g \pa_y \eh
\eee
yield the hoped bound for the linear term. We estimate afterwards the non linear term. We know that:
\bee
\hat N(\eh) = 3Q \eh^2 + 3 \hat f \eh^2 + \eh^3.
\eee 
where we note
\bee
 \hat f =3(\hat Q_b - Q)
\eee
we are within the bounds of the profiles $(T_i)_{1 \leq i \leq 3}$:
\bee
| \hat f + \pa_y \hat f + H \hat f  | &\lesssim& b \ \ \ \mbox{for} \ \ \ y \leq 1 \\
\left| \frac{\hat f}{1+|\log y|} + y\pa_y \hat f + y^2H \hat f  \right| &\lesssim& b \ \ \ \mbox{for} \ \ \ y \geq 1.
\eee
With the estimates $L^{\infty}$ of the Lemma \ref{lemmainterpolation}, of the second localisation \fref{radiation} and \fref{radiationbis}, we have:
\bee
\int |H \hat N(\eh)|^2 &\lesssim& \frac{\Delta(\eh^2)^2}{1+y^2} + \frac{(\pa_y\eh^2)^2}{1+y^6} + \int \frac{\eh^4}{1+y^8}\\
&+&\left\| \eh \right\|_{L^{\infty}}^2 \int |H \hat f \eh|^2 + |\Delta \eh \hat f |^2 + |\pa_y \hat f \pa_y \eh|^2 \\
&+& \left\| \eh \right\|_{L^{\infty}}^4 \hat{\mathcal E}_2 + \left\| \eh \right\|_{L^{\infty}}^2 \int |\Delta \eh^2|^2 + \left\| \eh \right\|_{L^{\infty}}^2\left\|\pa_y \eh \right\|_{L^{\infty}}^2 \int |\pa_y \eh|^2 \\
&\lesssim& b^5.
\eee
Proposition \ref{ttpptt} is proved.
\section{Proof of the proposition \ref{bootstrapregime}}
\subsection{Improved bound} The twice Lyapounov monotonicity properties give us the arguments to get better the priori bounds, under the a priori control \fref{init4h} on the unstable direction. 
\begin{lemma}[Improved bounds under the a priori control \fref{init4h}]
\label{improvedbounds}
Assume that K in \fref{init1h}, \fref{init2h}, \fref{init3h}, and \fref{init3hbis} has chosen large enough. Then, $\forall t \in [0,t_1]$:
\be
\label{improved1}
0 \leq b(t) \leq 2b(0),
\ee
\be
\label{improved2}
\int |\nabla \e|^2 \leq \sqrt{b(0)},
\ee
\be
\label{improved3}
|\mathcal E_2(t)| \leq \frac{K}{2}b^2(t)|\log b(t)|^5,
\ee
\be
\label{improved3bis}
|\mathcal E_4(t)| \leq \frac K2 \frac{b^4(t)}{|\log b(t)|^2}.
\ee
\end{lemma}
{\bf Proof of the Lemma \ref{improvedbounds}}\\ \par
{\bf Step 1} Positivity and smallness of $b(t)$
\\
\\
The proof of \fref{improved1} is a direct consequence of the modulation equations. Indeed, this last equation \fref{parameterspresicely}
 yields that
\be
b_s \leq 0
\ee
We must now prove that $b(t)$ can't be negative. We argue by contradiction. As $b(0) \geq 0$ and b(t) is a continue function, we suppose that it exists such $t_0$ that $b(t_0) = 0$. With the modulation equations, we have that:
\be
|b_s|\leq 2 b^2
\ee
Hence, there exists 
$\delta$ such that $b(t) = 0$ on $[t_0 - \delta, t_0]$, and thus from \fref{init3hbis}, $\lambda(t) = \lambda(t_0)$ and $u(t) = Q_{\lambda(t_0)}$ on $[t_0 - \delta, t_0]$. Iterating on $\delta >0$, we conclude that $u_0$ is initially a rescaling of Q, meaning a contradiction. 
\\
\par
{\bf Step 2} Energy bound
\\
\\
\fref{improved2} is a consequence of the decrease of energy. Indeed, let
\be
\tilde{\e} = \tilde{\alpha} + \e.
\ee
Then
\bea
\nonumber E(u)& =& \frac12 \int |\nabla u|^2 - \frac14 \int |u|^4\\
\nonumber&=& \frac12 \left\{\int |\nabla Q|^2 + \int |\nabla \tilde {\e}|^2\right\} + \int \pa_y Q\pa_y \tilde{\e} - \frac14\int \left[ Q^4 + 4 Q^3\tilde{\e} + 6 Q^2\tilde{\e}^2 + 4Q\tilde{\e}^3 + \tilde{\e}^4 \right]\\
\label{energybound1}&=& E(Q) + (H \tilde{\e},\tilde{\e}) - \frac14\int \left[ 4Q\tilde{\e}^3 + \tilde{\e}^4 \right].
\eea
Now,
\bee
(H \tilde{\e},\tilde{\e}) &=& (H\e, \e) + (H \tilde{\alpha},\tilde{\alpha}) + 2 (\tilde{\alpha}, H\e) \\
&=&(H\e,\e) + O(b|\log b|^C).
\eee
The last equality comes from Cauchy-Schwarz, the bound \fref{init3h} for $\e_2$ and the inequalities:
\bee
\left \| \tilde{\alpha} \right\|^2_{L^2} \lesssim |\log b|^C \ \ \ \left \| H\tilde{\alpha} \right\|^2_{L^2} \lesssim b^4 |\log b|^C.
\eee
Moreover, using \fref{coercivityHalone}, the orthogonality conditions \fref{ortho}, and the \begin{it}a priori\end{it} bound on the unstable direction \fref{init4h}
 \bee
 (H\e, \e) \geq c\int |\nabla \e|^2  - \frac 1c \left[ (\e, \Phi_M)^2 +(\e,\psi)^2 \right] \gtrsim \int |\nabla \e|^2 + O\left(\frac{b^5}{|\log b|^2}\right).
 \eee
 Thus,
\be
\label{energybound2}
(H \tilde{\e},\tilde{\e}) \gtrsim  \int |\nabla \e|^2 + O(b|\log b|^C).
\ee 
Let us see the nonlinear terms. We recall that:
\bee
\| \tilde{\e} \|_{L^{\infty}} \lesssim \| \e \|_{L^{\infty}} + \| \tilde{\alpha} \|_{L^{\infty}} \lesssim b|\log b| \ \ \ \left \|\frac{ \tilde{\alpha}}{y} \right\|^2_{L^2} + \left \| \nabla\tilde{\alpha} \right\|^2_{L^2} \lesssim b^2|\log b|^C.
\eee
Like this, with the bound \fref{init2h} and the fact that:
\be
\label{youpi}
\forall u \in H^1_{rad}(\mathbb R^4), \ \ \ \int \frac{|u|^2}{y^2} \lesssim \int |\nabla u|^2, \hspace{0,3cm}\left( \int |u|^4\right)^\frac{1}{4} \lesssim \left(\int |\nabla u|^2\right)^\frac12,
\ee
thus
\bea
\nonumber \int \left[ 4Q\tilde{\e}^3 + \tilde{\e}^4 \right] &\lesssim& \| \tilde{\e} \|_{L^{\infty}}  \left(\int \frac{|\e|^2}{1+y^2} + \left \| \frac{ \tilde{\alpha}}{y} \right\|^2_{L^2} \right) + \left(\int |\nabla \e|^2\right)^2 + \left(\int |\nabla \tilde{\alpha}|^2\right)^2\\
\label{energybound3}&\lesssim& \sqrt{b(0)} \int |\nabla \e|^2 + O(b |\log b|^C).
\eea
The first inequality of \fref{youpi} comes from the Lemma \fref{lemmetropcool}. The second is a classical result of Sobolev, Gagliardo and Nirenberg in dimension $N=4$. This result in general case and its proof is available in \cite{BH} with the Theorem IX.9. By construction, 
\be
\label{energybound4}
0 \leq E(u)-E(Q) \lesssim b^2(0)|\log b(0)|^C.
\ee
Injecting \fref{energybound2},\fref{energybound3} and \fref{energybound4} into \fref{energybound1} concludes the proof of \fref{improved2}.
\\
\par
{\bf Step 3} Control of $\mathcal E_4$ \\
\\
We argue similarily as in \cite{RS}. 
$\forall t\in [0,t_1)$,
\bea
\label{monotnyintegree}
\mathcal E_4(t) & \leq & 2\left(\frac{\lambda(t)}{\lambda(0)}\right)^6\left[\mathcal E_4(0)+C\sqrt{b(0)}\frac{b^4(0)}{|\log b(0)|^2}\right]+\frac{b^4(t)}{|\log b(t)|^2}\\
\nonumber & + & C\left[1+\frac{K}{\log M}+\sqrt{K}\right]\lambda^6(t)\int_0^t\frac{b}{\lambda^8}\frac{b^4}{|\log b|^2}d\tau
\eea
for some universal constant $C>0$ independent of $M$.\\
Let us now consider two constants 
\be
\label{conno}
\alpha_1=2-\frac{C_1}{\sqrt{\log M}}, \ \ \alpha_2=2+\frac{C_2}{\sqrt{\log M}}
\ee for some large enough universal constanst $C_1,C_2$. We compute using the modulation equations \fref{parameters}, \fref{parameterspresicely} and the bootstrap bound \fref{init3hbis}:
\bee
\frac{d}{ds}\left\{\frac{|\log b|^{\alpha_i}b}{\lambda}\right\}& = & \frac{|\log b|^{\alpha_i}}{\lambda}\left[\left(1-\frac{\alpha_i}{|\log b|}\right)b_s-\frac{\lambda_s}{\lambda}b\right]\\
& = & \frac{|\log b|^{\alpha_i}}{\lambda}\left[\left(1-\frac{\alpha_i}{|\log b|}\right)b_s+b^2 + O\left(\dfrac{b^3}{|\log b|}\right)\right]\\
& = & \left(1-\frac{\alpha_i}{|\log b|}\right)\frac{|\log b|^{\alpha_i}}{\lambda}\left[b_s+b^2\left(1+\frac{\alpha_i}{|\log b|}+O\left(\frac{1}{|\log b|^2}\right)\right)\right]\\
& &\left\{ \begin{array}{ll}\leq 0 \ \ \mbox{for} \ \ i=1\\\geq 0 \ \ \mbox{for}  \  \ i=2.
\end{array}\right .
\eee
Integrating this from $0$ to $t$ yields: 
\be
\label{lawintegrationone}
\frac{b(0)}{\lambda(0)}\left(\frac{|\log b(0||}{|\log b(t)|}\right)^{\alpha_2}
\leq \frac{b(t)}{\lambda(t)}\leq \frac{b(0)}{\lambda(0)}\left(\frac{|\log b(0||}{|\log b(t)|}\right)^{\alpha_1}.
\ee
This yields in particular using the initial bound \fref{init2} and the bound \fref{init1h}:
\be
\label{estationtermezero}
\left(\frac{\lambda(t)}{\lambda(0)}\right)^6\mathcal E_4(0)\leq (b(t)|\log b(t)|^{\alpha_2})^6\frac{\mathcal E_0}{(b(0)|\log b(0)|^{\alpha_2})^6}\leq \frac{b^4(t)}{|\log b(t)|^2},
\ee 
\bea
\label{estationtermezerobis}
\nonumber C\left(\frac{\lambda(t)}{\lambda(0)}\right)^6\sqrt{b(0)}\frac{b^4(0)}{|\log b(0)|^2} & \lesssim &\left(\frac{b(t)|\log b(t)|^{\alpha_2}}{b(0)|\log b(0)|^{\alpha_2}}\right)^6\sqrt{b(0)}\frac{b^4(0)}{|\log b(0)|^2}\\
& \lesssim & C(b(t))^{4+\frac14}\leq \frac{b^4(t)}{|\log b(t)|^2}.
\eea
We now compute explicitely using $b=-\lambda\lambda_{t}+O\left(\frac{b^2}{|\log b|}\right)$ from \fref{parameters}:
\bee
\int_0^t\frac{b}{\lambda^8}\frac{b^4}{|\log b|^2}d\sigma& = & \frac16\left[\frac{b^4}{\lambda^6|\log b|^2}\right]_0^t-\frac 16\int_0^{t}\frac{b_tb^3}{\lambda^6|\log b|^2}\left(4+\frac{2}{|\log b|}\right)d\tau\\
& + & O\left(\int_0^t\frac{b}{\lambda^8}\frac{b^5}{|\log b|^2}d\tau\right)
\eee
which implies using now $|b_s+b^2|\lesssim \frac{b^2}{|\log b|}$ from \fref{parameterspresicely} and \fref{init3hbis}:
$$\lambda^6(t)\int_0^t\frac{b}{\lambda^8}\frac{b^4}{|\log b|^2}d\sigma \lesssim \left[1+O\left(\frac{1}{|\log b_0|}\right)\right]\frac{b^4(t)}{|\log b(t)|^2}.$$ Injecting this together with \fref{estationtermezero}, \fref{estationtermezerobis} into \fref{monotnyintegree} yields $$\mathcal E_4(t)\leq C \frac{b^4(t)}{|\log b(t)|^2}\left[1+\frac{K}{\log M}+\sqrt{K}\right]$$ for some universal constant $C>0$ independent of $K$ and $M$, and thus \fref{improved3bis}  follows for $K$ large enough independent of $M$.\\
\par
{\bf Step 4} Control of $\mathcal E_2$\\
\\
Similarly to the control of $\mathcal E_4$, we give the same proof as well as in \cite{RS}. We integrate the monotonicity formula \fref{monotonicityE2} after recalling the estimate of the difference between $\hat {\mathcal E}_2$ and $\mathcal E_2$:
\bea
\label{cnikheoeif}
\mathcal E_2(t) & = &  \lambda^2(t)\|w_2(t)\|^2_{L^2}\lesssim  \|H\zeta_b(t)\|^2_{L^2}+\lambda^2(t)\|\wh_2(t)\|^2_{L^2}\\
\nonumber & \lesssim & b^4(t)|\log b(t)|^2+\left(\frac{\lambda(t)}{\lambda(0)}\right)^2\left[\mathcal E_2(0)+b^2(0)|\log b(0)|^2\right]+\lambda^2(t)\int_0^t\frac{b^3|\log b|^2}{\lambda^4(\tau)}d\tau.
\eea
From \fref{init2}, \fref{lawintegrationone}:
\bee
\left(\frac{\lambda(t)}{\lambda(0)}\right)^2\left[\mathcal E_2(0)+b^2(0)|\log b(0)|^2\right] & \lesssim & \frac{(b(0))^{10}+b^2(0)|\log b(0)|^2}{(b(0)|\log b(0)|^{\alpha_2}|)^2} b^2(t)|\log b(t)|^{2\alpha_2}\\
& \leq & b^2(t)|\log b(t)|^{4+\frac14},
\eee
We now use the bound $b_s\lesssim -b^2$ and \fref{lawintegrationone} to estimate: 
\bee
&&\lambda^2(t)\int_0^t\frac{b^3|\log b|^2}{\lambda^4(\tau)}d\tau\lesssim  \lambda^2(t)\int_0^t\frac{-b_tb|\log b|^2}{\lambda^2(\tau)}d\tau\\
& \lesssim & \left(\frac{\lambda(t)}{\lambda(0)}\right)^2b^2(0)|\log b(0)|^{2\alpha_1}\int_0^t\frac{-b_t}{b|\log b|^{2\alpha_1-2}}d\tau\\
& \lesssim &  \left(\frac{\lambda(t)}{\lambda(0)}\right)^2b^2(0)|\log b(0)|^{2\alpha_1}\frac{1}{|\log b(0)|^{2\alpha_1-3}}\\
& \lesssim & b^2(t)|\log b(t)|^{2\alpha_2}\frac{|\log b(0)|^3}{|\log b(0)|^{2\alpha_2}}\lesssim b^2(t)|\log b(t)|^{4+\frac 14}.
\eee
Injecting these bounds into \fref{cnikheoeif} yields: $$\mathcal E_2(t)\lesssim b^2(t)|\log b(t)|^{4+\frac14}$$ and concludes the proof of \fref{improved3}.\\

\subsection{Dynamic of the unstable mode}
To conclude the Proposition \ref{bootstrapregime}, we must study the dynamic of the unstable mode, which is the object of this subsection . We recall that $ \kappa (t) = (\varepsilon(t),\psi)$.
\begin{lemma}[Control of the unstable mode]
\label{dynamicunstablemode}
There holds: for all $t \in [0,T_1(a^+)],$
\be
\left| \frac{d\kappa}{ds} - \zeta \kappa \right| \leq \sqrt b \frac{b^{\frac52}}{|\log b|}.
\ee
\end{lemma}
{\bf Proof of the Lemma \ref{dynamicunstablemode}}
\\
\par
We compute the equation satisfied by $\kappa$ by taking the inner product of \fref{eqofepsilon} with the well localized direction $\psi$ to get:
\be
\label{eqofkappa}
 \frac{d\kappa}{ds} - \zeta \kappa = E(\varepsilon)
\ee
with
\be
\label{defEkappa}
E(\varepsilon) = (-\tilde{\Psi}_b,\psi) + (L(\varepsilon),\psi) + (N(\varepsilon),\psi) - (Mod,\psi) +\frac{\lambda_s}{\lambda} (\Lambda \varepsilon,\psi).
\ee
We now estimate all terms of RHS. We recall the exponential localization of $\psi$ as well as the orthogonality $(\psi,\Lambda Q)=0$. 
To begin, using \fref{controleh4erreurtilde}
\be
\label{psibkappa}
\left|(\tilde{\Psi}_b,\psi)\right| = \left|\frac{1}{\zeta^2}(\tilde{\Psi}_b,H^2\psi)\right| = \left|\frac{1}{\zeta^2}(H^2\tilde{\Psi}_b,\psi)\right| \lesssim \left(\int|H^2(\tilde{\Psi}_b)|^2\right)^{\frac{1}{2}} \lesssim  \sqrt b \frac{b^{\frac52}}{|\log b|}.
\ee
From the definition \fref{defLepsilon} of $L(\varepsilon)$, we have the following bound:
\bee
\left|L(\varepsilon) \right| \lesssim b y^{10} |\varepsilon|.
\eee
Thus,
\be
\left| (L(\varepsilon),\psi) \right| \lesssim b \left| \left(\frac{|\varepsilon|}{y^2(1+y^2)(1 + |\log y|)},(1+y^{14})(1+|\log y|)\psi\right) \right| \lesssim  \sqrt b \frac{b^{\frac52}}{|\log b|}.
\ee
In the same way, with \fref{defNepsilon},
\bee
\left|N(\varepsilon) \right| \lesssim \left(b y^{10} |\varepsilon| +  \left\| \frac{\varepsilon}{1+y}\right\|^2_{L^{\infty}} y^2 \right)|\varepsilon|.
\eee
\fref{normeinfty} and the fact that $\forall i, \left\| y^i \psi \right\|_{L^{\infty}} \lesssim 1 $ yield:
\be
\left| (N(\varepsilon),\psi) \right|  \lesssim  \sqrt b \frac{b^{\frac52}}{|\log b|}.
\ee
\\
With the notation of the Lemma \ref{localisation}, we have
\bee
\left| (Mod,\psi) \right| &=&  \left|\left(- \left( \frac{\lambda_s}{\lambda} +b\right) \Lambda \qbt + (b_s + b^2)(\tilde{T_1} +2b\tilde{T_2}),\psi\right) \right|\\
&\lesssim& \left|  \frac{\lambda_s}{\lambda} +b\right| \left| \left(\Lambda Q + \Lambda \tilde{\alpha},\psi\right)\right| + \left| b_s + b^2\right|\left|\left( \tilde{T_1} +2b\tilde{T_2},\psi\right)\right|.
\eee
But
\bee
\left| \Lambda \tilde{\alpha}\right| + \left| 2b\tilde{T_2}\right| \lesssim by^{10}
\eee
and
\bee
\left|\left( \tilde{T_1},\psi\right)\right| = \left|\frac{-1}{\zeta}\left( \tilde{T_1},H\psi\right)\right| = \left|\frac{-1}{\zeta}\left( H\tilde{T_1},\psi\right)\right| \lesssim \left|\int_{y \geq B_1}\left(H\tilde{T_1}- \Lambda Q\right)\psi \right| \lesssim b.
\eee
Hence, with the modulation equations
\be
\label{modkappa}
\left| (Mod,\psi) \right| \lesssim  \sqrt b \frac{b^{\frac52}}{|\log b|}.
\ee
For the last term, we use \fref{estun}.
\be
\left| \frac{\lambda_s}{\lambda} (\Lambda \varepsilon,\psi)\right| = \left| \frac{\lambda_s}{\lambda}\right|\left| \left(\frac{\Lambda \varepsilon}{y^2(1+y^2)(1+|\log y|)},y^2(1+y^2)(1+|\log y|)\psi\right)\right|\lesssim  \sqrt b \frac{b^{\frac52}}{|\log b|}.
\ee
This concludes the proof of the Lemma \ref{dynamicunstablemode}.
\subsection{Conclusion}
We have at our disposal all the elements to finish the proof using a rough argument which does not give a sharp information on the link between the initial data and the choice of $a^+$, and uniqueness for example is not covered at this stage. The keystone of the proof is the fact that the map:
\bee
\left[ -2\frac{b_0^{\frac 52}}{|\log b_0|} ; 2\frac{b_0^{\frac 52}}{|\log b_0|}\right]&\rightarrow&\mathbb R^+\\
a^+&\rightarrow & T_1(a^+)
\eee
is continuous as a consequence of the strictly outgoing behavior on exit \fref{eqofkappa} defined by \fref{init4h}. This classical argument is displayed in details  in \cite{CMM},  Lemma 6, in a more complicated setting and therefore left to the reader. Hence, we also have the continuity of the map 
\bee
\left[ -2\frac{b_0^{\frac 52}}{|\log b_0|} ; 2\frac{b_0^{\frac 52}}{|\log b_0|}\right]&\rightarrow&\mathbb R^+\\
a^+&\rightarrow & \kappa (T_1(a^+)).
\eee
In agreement with the dynamics of the unstable mode found in the last subsection, we know that, for $a^+ = 2\frac{b_0^{\frac 52}}{|\log b_0|} $:
\bee
\frac{d}{ds}\kappa(0) = 2\zeta \frac{b_0^{\frac 52}}{|\log b_0|} + O\left(\frac{b_0^3}{|\log b_0|}\right) > 0
\eee 
and then
\be
\kappa\left(T_1\left(2\frac{b_0^{\frac 52}}{|\log b_0|}\right)\right) = 2\frac{b_0^{\frac 52}}{|\log b_0|}.
\ee
Likewise,
\be
\kappa\left(T_1\left(-2\frac{b_0^{\frac 52}}{|\log b_0|}\right)\right) = -2\frac{b_0^{\frac 52}}{|\log b_0|}.
\ee
By contuinity, it exists $a^+ \in \left] -2\frac{b_0^{\frac 52}}{|\log b_0|} ; 2\frac{b_0^{\frac 52}}{|\log b_0|}\right[$ such that 
\be
\label{conditionsinequanon}
\kappa\left(T_1\left(a^+\right)\right)=0.
\ee
In addition, according to the definition of exit time and the Lemma \ref{improvedbounds}, we have two choices.  Either $|\kappa(T_1(a^+)) | = 2\frac{b(T_1(a^+))^{\frac 52}}{|\log b(T_1(a^+))|}$ or $T_1(a^+)$ is the life time of the solution.  If the first possibility is the good one, the condition \fref{conditionsinequanon} gives  $$2\frac{b(T_1(a^+))^{\frac 52}}{|\log b(T_1(a^+))|} =0.$$ As we have proved that $b(t)>0$ for $t < T$, where T is the life time of the solution, we have thus the second possibility. Notice that in this case, the two choices tally. This is exactly the Proposition \ref{bootstrapregime}. \\


\begin{appendix}

\section{$L^2$ coercivity estimates}
In this appendix, we let prove at the first Hardy inequalities for functions $u \in H^2_{rad}(\mathbb R^4)$. We will use afterwards this results to establish properties of weighted sub-coercivity for $H$ and $H^2$, which allow us to obtain coercive estimates for these operators under additional orthogonality conditions. This coercive estimates are crucial in our study. The proof lies in the continuation of the analysis in \cite{HR}.
\subsection{Hardy inequalities}
\begin{lemma}
\label{lemmetropcool}
There exists a constant C for which there holds, for any $v \in H^1_{rad}(\mathbb R^4)$ 
\be
\label{lemme1}
\left[ \int_{\mathbb R^4} \frac{v(y)^2}{y^2}\right]^{\frac12} + \underset{y \in \mathbb R^4}{\sup} \left( |yv(y)|\right) \leq C \left[ \int_{\mathbb R^4} |\nabla v(y)| ^2 \right]^{\frac12},
\ee
\be
\label{normeinftysup1}
 \left \| v \right\|^2_{L^{\infty}_{y \geq 1}} \lesssim \int_{y \geq 1}\left( \frac{|\nabla v|^2}{y^2} + \frac{|v|^2}{y^4}\right) + \int_{\frac12 \leq y \leq 1} \frac{|v|^2}{y^2}
\ee

\end{lemma}
{\it{Proof of the Lemma \ref{lemmetropcool}:}} 
By integration-by-parts:
\bee
 \int \frac{v(y)^2}{y^2}& =& \left[ \frac{v(y)^2y^2}2\right]^{\infty}_0 -  \int \frac{v(y)\pa_yv(y)}{y}\\
 & \lesssim & \left[ \int \frac{v(y)^2}{y^2}\right]^{\frac12}\left[ \int|\nabla v(y)| ^2 \right]^{\frac12}.
\eee
Next:
$$|v(y)|\lesssim \int_y^{+\infty}\frac{|\pa_yv|}{y^3}\lesssim \frac{1}y \left(\int|\pa_yv|^2\right)^{\frac 12},$$
\bee
|v^2(y)|\lesssim \int_{1\leq y\leq 2}|v|^2+\int_{y\geq 1}\frac{|v||\pa_yv|}{y^3}\lesssim \left(\int_{y\geq 1}\frac{|\pa_yv|^2}{y^2}\right)^{\frac 12}\left(\int \frac{|v|^2}{y^4}\right)^{\frac 12}.
\eee
This concludes the proof of Lemma \ref{lemmetropcool}. \\
\begin{lemma}[Hardy inegalities]
\label{Hardylemma}
$\forall R > 2$, $\forall v \in H^2_{rad}(\mathbb R^4)$, $\forall \gamma > 0$ there holds the following controls:
\bea
\label{hardy1}\int |\pa_{yy} v|^2 + \int \frac{|\pa_y v| ^2}{y^2} &\lesssim& \int (\Delta v)^2, \\
\label{hardy2}\int_{y \leq R} \frac{|v|^2}{y^4(1+|\log y|)^2} &\lesssim&  \int_{y \leq R} \frac{|\pa_y v| ^2}{y^2} +  \int_{1 \leq y \leq 2}|v|^2,\\
\label{hardy3}\int_{1 \leq y \leq R} \frac{|v|^2}{y^{4+\gamma}(1+|\log y|)^2} &\lesssim&  \int_{1 \leq y \leq R} \frac{|\nabla v| ^2}{y^{2+\gamma}(1+|\log y|)^2} + C_{\gamma} \int_{1 \leq y \leq 2}|v|^2.
\eea
\end{lemma}
{\it Proof of the Lemma \ref{Hardylemma}:}. Let $v$ smooth and radially symmetric. \fref{hardy1} follows from the explicit formula after integration of parts
\bee
\int (\Delta v)^2 = \int \left( \pa_{yy} v + \frac3y \pa_y v \right)^2 = \int |\pa_{yy} v|^2 + 3 \int \frac{|\pa_y v| ^2}{y^2}.
\eee
To prove \fref{hardy2} and \fref{hardy3}, from the one dimensional Sobolev embedding $H^1(1 \leq y \leq 2)$ in $L^{\infty}(1 \leq y \leq 2)$, we obtain
\be
\label{inegalitev1}
|v(1)|^2 \lesssim \int_{1 \leq y \leq 2} \left( |v|^2 + |\pa_y v|^2 \right).
\ee
Let $f(y) = -\frac{{\bf e}_y}{y^3(1+\log y)}$ so that $\nabla . f = \frac{1}{y^4(1+\log y)^2}$, and integrate by parts to get:
\bea
\nonumber &&\int_{1 \leq y \leq R} \frac{|v|^2}{y^4(1+|\log y|)^2} = \int_{1 \leq y \leq R} |v|^2 \nabla .f  \\
\nonumber &=& - \left[ \frac{|v|^2}{1+\log y}\right]^R_1 + 2 \int_{y \leq R} \frac{v \pa_y v}{y^3(1+\log y)} \\
\label{zz}&\lesssim& |v(1)|^2 + \left(\int_{y \leq R} \frac{|v|^2}{y^4(1+|\log y|)^2}\right)^{\frac12} \left(\int_{y\leq R} \frac{|\pa_y v| ^2}{y^2}\right)^{\frac12}.
\eea
Similarly, using $\tilde f(y) = -\frac{{\bf e}_y}{y^3(1-\log y)}$, we get:
\bea
\nonumber &&\int_{\varepsilon \leq y \leq 1} \frac{|v|^2}{y^4(1- \log y)^2} = \int_{\varepsilon \leq y \leq 1} |v|^2 \nabla .f  \\
\nonumber &=& \left[ \frac{|v|^2}{1-\log y}\right]^1_{\varepsilon} + 2 \int_{y \leq 1} \frac{v \pa_y v}{y^3(1-\log y)} \\
\label{zzz}&\lesssim& |v(1)|^2 + \left(\int_{y \leq R} \frac{|v|^2}{y^4(1+|\log y|)^2}\right)^{\frac12} \left(\int_{y\leq R} \frac{|\pa_y v| ^2}{y^2}\right)^{\frac12}.
\eea
\fref{inegalitev1}, \fref{zz} and \fref{zzz} now yield \fref{hardy2}. To prove \fref{hardy3}, let $\gamma > 0$, and 
\bee
f(y) = -\frac{{\bf e}_y}{y^{\gamma+3}(1+\log y)^2}
\eee
so that for $y \geq 1$
\bee
\nabla.f(y) = \frac{1}{y^{\gamma+4}(1+\log y)^2}\left[ \gamma + \frac2{1+\log y}\right] \geq \frac{\gamma}{y^{\gamma+4}(1+\log y)^2}.
\eee
We then integrate by parts to get:
\bee
&&\gamma \int_{1 \leq y \leq R} \frac{|v|^2}{y^{4+\gamma}(1+|\log y|)^2} \leq \int_{1 \leq y \leq R} |v|^2 \nabla.f\\
&\leq& - \left[ \frac{|v|^2}{y^{3+\gamma}(1+|\log y|)^2}\right]^R_1 + 2 \int_{1\leq y \leq R} \frac{|v \pa_y v|}{y^{\gamma+3}(1+\log y)^2}\\
&\leq& C \int_{1\leq y \leq 2} |v|^2 + 2 \left(\int_{ y \leq R} \frac{|v|^2}{y^{4+\gamma}(1+|\log y|)^2}\right)^{\frac{1}{2}}\left( \int_{y \leq R} \frac{|\nabla v| ^2}{y^{2+\gamma}(1+|\log y|)^2}\right)^{\frac{1}{2}}.
\eee
and \fref{hardy3} follows. 
\subsection{Sub-positivy estimates with H}
The following lemma highlights the negative part of the operator $H$. We recall that this operator possesses a unique nonpositive direction $\psi$.
\begin{lemma}
\label{subpositivity}
Let $u \in H^2_{rad}(\mathbb R^4)$, then there exists a constant $C>0$ such that:
\be
\left(Hu,u\right) \geq -C\left(u,\psi \right)^2.
\ee
\end{lemma}
{\it Proof of the Lemma \ref{subpositivity}:} Let $u \in H^2_{rad}(\mathbb R^4)$. There exists an unique decomposition of $u$:
\bee
u = \kappa \psi + v
\eee
with the orthogonality condition
\bee
(\psi,v)=0.
\eee
By definition, we have
\bee
\kappa = \frac{\left(u,\psi\right)}{\left(\psi,\psi\right)}.
\eee
Morever, the uniqueness of the negative direction of H gives
\bee
\left(Hv,v\right) \geq 0.
\eee
Thus,
\bee
\left(Hu,u\right) &=& \kappa^2(H\psi,\psi) + \kappa\left[(Hv,\psi) + (v,H\psi)\right] + (Hv,v) \\
&=& \kappa^2(H\psi,\psi)  - 2\zeta \kappa (v,\psi) + (Hv,v) \\
&\geq& \frac{(H\psi,\psi)}{\left(\psi,\psi\right)^2}\left(u,\psi\right)^2 \\
&\geq& \frac{-\zeta}{\left(\psi,\psi\right)}\left(u,\psi\right)^2.
\eee
\subsection{Sub-coercivity estimates}
In this subsection, we prove sub-coercivity estimates for $H$ and $H^2$ which are the key to the proof of coercive estimates for these operators under additional orthogonality conditions.
\begin{lemma}[Sub-coercivity estimates with H]
\label{lemme1forH}
Let $u \in H^2_{rad}(\mathbb R^4)$, then there exists  constants $\delta>0, C>0$ such that:
\be
\label{coercivityHalone}
(H\e,\e) \geq c \int |\nabla \e|^2 - \frac 1c \left[ (\e, \Phi_M)^2 + (\e, \psi)^2\right],
\ee
\be
\label{propHcoercivity}
\int |\pa_{yy} u|^2 + \int \frac{|\pa_y u|^2}{y^2} + \int \frac{u^2}{y^4(1+|\log y|)^2} - C\left[\int \frac{|\pa_y u|^2}{1+y^4} + \int \frac{u^2}{1+y^8}\right] \lesssim \int (Hu)^2.
\ee
\end{lemma}
{\it Proof of the Lemma \ref{lemme1forH}:} \fref{propHcoercivity} is a direct consequence of the inegalities \fref{hardy1}, \fref{hardy2} and the following decomposition:
\bee
\int (Hu)^2 &=& \int(\Delta u + Vu)^2 = \int (\Delta u)^2 - 2 \int V(\pa_yu)^2+ \int (\Delta V + V^2)u^2\\
&\gtrsim& \left[ \int (\Delta u)^2 + \int \frac{u^2}{1+y^6}\right] - C \left[\int \frac{|\pa_y u|^2}{1+y^4} + \int \frac{u^2}{1+y^8}\ \right]
\eee
where we used that
\bee
V(y) = \frac{192}{1+y^4} \left[ 1 + O\left(\frac1{1+y^2}\right)\right] \ \ \mbox{as} \ \ y\to+\infty.
\eee
\begin{lemma}[Weighted sub-coercivity for H]
\label{lemme2forH}
Let $u \in  H^4_{rad}(\mathbb R^4)$, then exists a constant C such that:
\bea
\nonumber &&\int \frac{|u|^2}{y^4(1+y^4)(1+|\log y|)^2} + \int \frac{|\pa_y u|^2}{y^6(1+|\log y|)^2} + \int \frac{|\pa_{yy} u|^2}{y^4(1+|\log y|)^2}  \\
\nonumber &+& \int \frac{|\pa_y^3 u|^2}{y^2(1+|\log y|)^2} +\int\frac{ |\pa_y^4 u|^2}{(1+|\log y|)^2} \\
\nonumber &-& C\left[ \int \frac{|u|^2}{y^2(1+y^8)(1+|\log y|)^2} + \int \frac{|\pa_y u|^2}{y^4(1+y^4)(1+|\log y|)^2}\right]\\
\label{reslemme2forH}&\lesssim&  \int \frac{|H u|^2}{y^4(1+|\log y|)^2} +\int \frac{|\pa_y Hu|^2}{y^2(1+|\log y|)^2} + \int \frac{|\pa_{yy}H u|^2}{(1+|\log y|)^2}.
\eea
\end{lemma}
\begin{remark} Using \fref{lemme1}, \fref{hardy2}, and \fref{hardy3}, 
$u \in  H^4_{rad}(\mathbb R^4)$ yields 
\bee
&&\int \frac{|u|^2}{y^4(1+y^4)(1+|\log y|)^2} + \int \frac{|\pa_y u|^2}{y^6(1+|\log y|)^2} + \int \frac{|\pa_{yy} u|^2}{y^4(1+|\log y|)^2} \\
&+&   \int \frac{|\pa_y^3 u|^2}{y^2(1+|\log y|)^2} +\int\frac{ |\pa_y^4 u|^2}{(1+|\log y|)^2} \\
&+& \int \frac{|H u|^2}{y^4(1+|\log y|)^2} +\int \frac{|\pa_y Hu|^2}{y^2(1+|\log y|)^2} + \int \frac{|\pa_{yy}H u|^2}{(1+|\log y|)^2}< \infty.
\eee
\end{remark}
{\it Proof of Lemma \ref{lemme2forH}}:
Let $\chi(y)$ be a smooth cut-off function with support in $y \geq 1$ and equal to 1 for $y \geq 2$.
\bee
&&\int \chi \frac{|H u|^2}{y^4(1+|\log y|)^2} = \int \chi \frac{\left[- \pa_y\left(y^3\pa_y u\right)  - y^3Vu\right]^2}{y^{10}(1+|\log y|)^2}\\
&=& \int  \chi \frac{ \left|\pa_y\left(y^3\pa_y u\right)\right|^2}{y^{10}(1+|\log y|)^2} + 2\int  \chi\frac{\pa_y\left(y^3\pa_y u\right) Vu}{y^{7}(1+|\log y|)^2} + \int  \chi\frac{V^2u^2}{y^{4}(1+|\log y|)^2} \\
&=& \int  \chi\frac{ \left|\pa_y\left(y^3\pa_y u\right)\right|^2}{y^{10}(1+|\log y|)^2} - 2\int \chi \frac{V(\pa_yu)^2 }{y^{4}(1+|\log y|)^2} + \int \chi \frac{V^2u^2}{y^{4}(1+|\log y|)^2} \\
&+& \int |u|^2  \Delta \left( \chi \frac V{y^4(1+|\log y|)^2}\right).
\eee
We now observe that for $k \geq 0$
\bee
|\pa_y^k V(y)| \lesssim \frac1{1+y^{4+k}}
\eee
and thus,
\bee
&&\left|- 2 \int \chi \frac{V(\pa_yu)^2 }{y^{4}(1+|\log y|)^2} + \int \chi \frac{V^2u^2}{y^{4}(1+|\log y|)^2} + \int  \chi u^2 \Delta\left( \frac V{y^4(1+|\log y|)^2}\right) \right| \\
&\lesssim& \int \frac{|u|^2}{y^2(1+y^8)(1+|\log y|)^2} + \int \frac{|\pa_y u|^2}{y^4(1+y^4)(1+|\log y|)^2}.
\eee
Hence,
\bee
\int \chi \frac{|H u|^2}{y^4(1+|\log y|)^2} &\gtrsim& \int \chi \frac{ \pa_y\left(y^3\pa_y u\right)^2}{y^{10}(1+|\log y|)^2} \\
&-& C \left[ \int \frac{|u|^2}{y^2(1+y^8)(1+|\log y|)^2} + \int \frac{|\pa_y u|^2}{y^4(1+y^4)(1+|\log y|)^2}\right].
\eee
We may apply twice the Hardy inegality \fref{hardy3} with $\gamma = 8$ and $\gamma = 4$ and get for a sufficiently large universal constant $R$:
\bee
 \int \chi \frac{ \pa_y\left(y^3\pa_y u\right)^2}{y^{10}(1+|\log y|)^2} &\gtrsim&  \int _{y \geq R}\frac{| \pa_y u|^2}{y^{6}(1+|\log y|)^2} - C\int \frac{|\pa_y u|^2}{y^4(1+y^4)(1+|\log y|)^2} \\
 &\gtrsim& \int_{y \geq R} \frac{|u|^2}{y^4(1+y^4)(1+|\log y|)^2} \\
 &-& C \left[ \int \frac{|u|^2}{y^2(1+y^8)(1+|\log y|)^2} + \int \frac{|\pa_y u|^2}{y^4(1+y^4)(1+|\log y|)^2}\right].
\eee
and finally,
\bea
\nonumber &&\int _{y \geq 2}\frac{|u|^2}{y^4(1+y^4)(1+|\log y|)^2} + \int_{y \geq 2} \frac{|\pa_y u|^2}{y^6(1+|\log y|)^2} + \int_{y \geq 2} \frac{|\pa_{yy} u|^2}{y^4(1+|\log y|)^2}  \\
\nonumber &-& C\left[ \int \frac{|u|^2}{y^2(1+y^6)(1+|\log y|)^2} + \int \frac{|\pa_y u|^2}{y^4(1+y^4)(1+|\log y|)^2}\right]\\
\label{pp}&\lesssim&  \int \chi \frac{|H u|^2}{y^4(1+|\log y|)^2}.
\eea
Now, the control of the third derivate for $y \leq 1$ follows from:
\bea
\nonumber  \int \chi\frac{|\pa_y Hu|^2}{y^2(1+|\log y|)^2}\geq \int\chi \frac{|\pa_y (\Delta u +Vu)|^2}{y^2(1+|\log y|)^2} \gtrsim \int\chi \frac{|\pa_y^3 u|^2}{y^2(1+|\log y|)^2}&& \\
\label{ppp} - C \left[ \int\chi \frac{|u|^2}{(1+y^{12})(1+|\log y|)^2} + \int \chi\frac{|\pa_y u|^2}{y^6(1+|\log y|)^2} + \int \chi\frac{|\pa_{yy} u|^2}{y^4(1+|\log y|)^2} \right] &&
\eea
and the fourth derivate from:
\bea
\nonumber && \int \chi\frac{|\pa_{yy }Hu|^2}{(1+|\log y|)^2}\geq \int \chi\frac{|\pa_{yy} (\Delta u +Vu)|^2}{(1+|\log y|)^2}\\
 &\gtrsim& \int \chi\frac{|\pa_y ^4 u|^2}{(1+|\log y|)^2} - C \left[ \int  \chi\frac{|\pa_y^3 u|^2}{y^2(1+|\log y|)^2} + \int \chi\frac{|\pa_{yy} u|^2}{y^4(1+|\log y|)^2}\right]\\
\label{pppp}&-&C\left[\int \chi\frac{|u|^2}{(1+y^{12})(1+|\log y|)^2} + \int \chi\frac{|\pa_y u|^2}{y^6(1+|\log y|)^2}  \right].
\eea
\fref{pp},\fref{ppp} and \fref{pppp} yield \fref{reslemme2forH}, away from the origin. Let us study this control near the origin. Let $\zeta = (1-\chi)^{\frac12}$. With the Lemma \ref{lemme1forH}, we have that:
\bee
&\int& \zeta^2 \frac{|Hu|^2}{y^4(1+|\log y|^2)} \gtrsim \int \zeta^2|Hu|^2 \gtrsim |H\zeta u|^2 - C \int _{1 \leq y \leq 2} \left( |\pa_y u|^2 + |u|^2\right) \\
&\gtrsim& \int \frac{|\zeta u|^2}{y^4(1+|\log y|^2)} - C \int \frac{|\zeta u|^2}{y^2} \gtrsim \int_{y \leq 1} \frac{|u|^2}{y^4(1+|\log y|^2)} - C \int_{y \leq 2} \frac{| u|^2}{y^2}.
\eee
Now, by definition, we have:
\bee
Hu = - \frac 1{y^3}\frac{\pa }{\pa_y}\left( y^3 \pa_y u\right) - Vu.
\eee
Hence
\be
\pa_y u = -\frac1{y^3} \int_{\tau\leq y}(Vu + Hu).
 \ee
 We then estimate from Cauchy-Schwarz and Fubini:
 \bea
\nonumber  \int_{y \leq 1} \frac{|\pa_y u|^2}{y^6(1+|\log y|^2)} &=& \int_{y \leq 1 } \frac{|\pa_y u|^2}{y^{12}(1+|\log y|^2)} \left( \int _0^y Vu + Hu \right)^2 \\
\nonumber  &\lesssim& \int_{y \leq 1 } \frac{|\pa_y u|^2}{y^{2}(1+|\log y|^2)} \left( \int _0^y \frac{|Vu(\tau)|^2 + |Hu(\tau)|^2}{\tau^3} \right) \\
 \nonumber &\lesssim& \int _{\tau \leq 1}\frac{|Vu(\tau)|^2 + |Hu(\tau)|^2}{\tau^3} \left(\int_{\tau \leq y \leq 1}\frac{1}{y^{2}(1+|\log y|^2)} \right) \\
 &\lesssim& \int _{y \leq 1} \frac{|H u|^2}{y^4(1+|\log y|)^2}.
 \eea
Finally, for the control of the other derivates near the origin:
\bee
\pa_{yy} u = - Hu - 3 \frac{\pa_y u}{y} - Vu.
\eee
Thus,
\be
\int_{y \leq 1} \frac{|\pa_{yy} u|^2}{y^4(1+|\log y|^2)} \lesssim \int_{y \leq 1} \frac{|H u|^2}{y^4(1+|\log y|^2)} + \int_{y \leq 1} \frac{|\pa_{y} u|^2}{y^6(1+|\log y|^2)} + \int_{y \leq 1} \frac{|u|^2}{y^4(1+|\log y|^2)},
 \ee
\bea
\nonumber \int_{y \leq 1} \frac{|\pa_{y}^3 u|^2}{y^2(1+|\log y|^2)} &\lesssim& \int_{y \leq 1} \frac{|\pa_y H u|^2}{y^2(1+|\log y|^2)} + \int_{y \leq 1} \frac{|\pa_{yy} u|^2}{y^4(1+|\log y|^2)} \\
&+& \int_{y \leq 1} \frac{|\pa_y u|^2}{y^6(1+|\log y|^2)} + \int_{y \leq 1} \frac{|u|^2}{y^4(1+|\log y|^2)},
 \eea 
\bea
\nonumber \int_{y \leq 1} \frac{|\pa_{y}^4 u|^2}{(1+|\log y|^2)} &\lesssim& \int_{y \leq 1} \frac{|\pa_{yy} H u|^2}{(1+|\log y|^2)} + \int_{y \leq 1} \frac{|\pa_{y}^3 u|^2}{y^4(1+|\log y|^2)} + \int_{y \leq 1} \frac{|\pa_{yy} u|^2}{y^4(1+|\log y|^2)} \\
&+& \int_{y \leq 1} \frac{|\pa_y u|^2}{y^6(1+|\log y|^2)} + \int_{y \leq 1} \frac{|u|^2}{y^4(1+|\log y|^2)}.
 \eea 
This concludes the proof.
 \par
We now combine the results of Lemma \ref{lemme1forH} and the Lemma \ref{lemme2forH}
\begin{lemma}[Sub-coercivity for $H^2$] 
\label{wxwx}
Let $u \in  H^4_{rad}(\mathbb R^4)$. Then,
\bea 
\label{dsds}\int |H^2u|^2 &\gtrsim& \int \frac{|H u|^2}{y^4(1+|\log y|)^2} +\int \frac{|\pa_y Hu|^2}{y^2(1+|\log y|)^2} + \int |\frac{\pa_{yy}H u}{(1+|\log y|)^2}|^2 \\
\nonumber &+&\int \frac{|u|^2}{y^4(1+y^4)(1+|\log y|)^2} + \int \frac{|\pa_y u|^2}{y^6(1+|\log y|)^2} + \int \frac{|\pa_{yy} u|^2}{y^4(1+|\log y|)^2} \\ 
\nonumber &+& \int \frac{|\pa_y^3 u|^2}{y^2(1+|\log y|)^2} \int\frac{ |\pa_y^4 u|^2}{(1+|\log y|)^2}  -C  \int \frac{|u|^2}{y^2(1+y^8)(1+|\log y|)^2} \\
\nonumber&-& C \left[ \int \frac{|\pa_y u|^2}{y^4(1+y^4)(1+|\log y|)^2} +\int \frac{|\pa_y Hu|^2}{1+y^4} + \int \frac{Hu^2}{1+y^8} \right].
\eea
\end{lemma}
\subsection{Coercivity of $H^2$} We are now in position to derive the fundamental coercivity property of $H^2$ at the heart of our analysis.
\begin{lemma}[Coercivity of $H^2$]
\label{coercivityH2}
Let $M\geq 1$ be a large enough universal constant. Let $\Phi_M$ be given by \fref{defdirection}. Then there exists a universal constant $C(M)>0$ such that for all  $u \in  H^4_{rad}(\mathbb R^4)$ satisfiying the orthogonality conditions:
\bee
(u,\Phi_M) = 0, \ \ \ \ (Hu, \Phi_M)=0
\eee
there holds:
\bea
\nonumber &&\int \frac{|H u|^2}{y^4(1+|\log y|)^2} +\int \frac{|\pa_y Hu|^2}{y^2(1+|\log y|)^2} + \int \frac{|\pa_{yy}H u|^2}{(1+|\log y|)^2} \\ 
\nonumber &+& \int \frac{|u|^2}{y^4(1+y^4)(1+|\log y|)^2} \int \frac{|\pa_y u|^2}{y^6(1+|\log y|)^2} + \int \frac{|\pa_{yy} u|^2}{y^4(1+|\log y|)^2}  \\
&+& \int \frac{|\pa_y^3 u|^2}{y^2(1+|\log y|)^2} +  \int\frac{ |\pa_y^4 u|^2}{(1+|\log y|)^2} \leq C(M) \int |H^2(u)|^2.
\eea
\end{lemma}
{\it Proof of the Lemma \ref{coercivityH2}:} We argue by contradiction. Let $M>0$ fixed and consider a normalized sequence $u_n$ 
\bea
\nonumber &&\int \frac{|H u_n|^2}{y^4(1+|\log y|)^2} +\int \frac{|\pa_y Hu_n|^2}{y^2(1+|\log y|)^2} + \int \frac{|\pa_{yy}H u_n|^2}{(1+|\log y|)^2} \\ 
\nonumber &+& \int \frac{|u_n|^2}{y^4(1+y^4)(1+|\log y|)^2} + \int \frac{|\pa_y u_n|^2}{y^6(1+|\log y|)^2} + \int \frac{|\pa_{yy} u_n|^2}{y^4(1+|\log y|)^2}  \\
&+& \int \frac{|\pa_y^3 u_n|^2}{y^2(1+|\log y|)^2} +  \int\frac{ |\pa_y^4 u_n|^2}{(1+|\log y|)^2} =1
\eea
satisfying the orthogonality conditions
\bee
(u_n,\Phi_M) = 0, \ \ \ \ (Hu_n, \Phi_M)=0
\eee 
and
\be
\label{cond2}
 \int |H^2(u_n)|^2 \leq  \frac1n.
\ee
The normalization condition implies that the sequence $u_n$ is uniformly bounded in $H^4_{loc}$. As a consequence, we can assume that $u_n$ weakly converges in $H^4_{loc}$ to $u_{\infty}$. Morever, $u_{\infty}$ satisfies the equation
\bee
H^2u_{\infty} = 0\ \ \mbox{for}\ \ r\>0.
\eee
Integrating this ODE leads to 
\bee
Hu_{\infty} = \alpha \Lambda Q + \beta \G \ \ \mbox{for} \ \ r>0.
\eee
Using the condition $u_{\infty} \in H^4_{loc}$, we can determine that $\beta = 0$. Hence, the function $u_{\infty}$ can be written in the form
\bee
u_{\infty} = - \alpha T_1 + \gamma \Lambda Q + \delta \G.
\eee
The condition $u_{\infty} \in H^4_{loc}$ yields that $\delta = 0$. Passing trough the limit in the orthogonality conditions, using that $u_n$ converges to $u_{\infty}$ weakly in $ H^4_{loc}$, we conclude that $u_{\infty}$ satisfies 
\bee
(u_{\infty} ,\Phi_M) = 0, \ \ \ \ (Hu_{\infty} , \Phi_M)=0.
\eee 
We may therefore determine the constant $\alpha$ and $\gamma$ using \fref{defdirection}, \fref{ortho} which yield $\alpha = \gamma = 0$ and thus $u_{\infty}=0$.
\par
The sub-coercitivity bound \fref{dsds} together with \fref{cond2} ensures:
\bee
&& \frac 1n \gtrsim  \int |H^2(u_n)|^2 \gtrsim \int \frac{|H u_n|^2}{y^4(1+|\log y|)^2} +\int \frac{|\pa_y Hu_n|^2}{y^2(1+|\log y|)^2} + \int \frac{|\pa_{yy}H u_n|^2}{(1+|\log y|)^2} \\ 
 &+&  \int \frac{|u_n|^2}{y^4(1+y^4)(1+|\log y|)^2} + \int \frac{|\pa_y u_n|^2}{y^6(1+|\log y|)^2} + \int \frac{|\pa_{yy} u_n|^2}{y^4(1+|\log y|)^2}\\
\nonumber &+& \int \frac{|\pa_y^3 u_n|^2}{y^2(1+|\log y|)^2}  +\int\frac{ |\pa_y^4 u_n|^2}{(1+|\log y|)^2} -C\left[\int \frac{|\pa_y Hu_n|^2}{1+y^4} + \int \frac{Hu_n^2}{1+y^8} \right]\\
\nonumber&-&C \left[ \int \frac{|u_n|^2}{y^2(1+y^8)(1+|\log y|)^2} + \int \frac{|\pa_y u_n|^2}{y^4(1+y^4)(1+|\log y|)^2}\right].
\eee
Coupling this with the normalization condition we obtain that
\bee
\int \frac{|\pa_y Hu_n|^2}{1+y^4} + \int \frac{Hu_n^2}{1+y^8} + \int \frac{|u_n|^2}{y^2(1+y^8)(1+|\log y|)^2} + \int \frac{|\pa_y u_n|^2}{y^4(1+y^4)(1+|\log y|)^2} \geq c
\eee
for some positive constant $c>0$. Since $u_n$ weakly converges to $u_{\infty}$ in $H^4_{loc}$ on any compact subinterval of $y \in (0,\infty)$, we can pass to the limit to conclude
\bee
\int \frac{|\pa_y Hu_{\infty}|^2}{1+y^4} + \int \frac{Hu_{\infty}^2}{1+y^8} + \int \frac{|u_{\infty}|^2}{y^2(1+y^8)(1+|\log y|)^2} + \int \frac{|\pa_y u_{\infty}|^2}{y^4(1+y^4)(1+|\log y|)^2} \geq c.
\eee
This contradicts the established identity $u_{\infty} = 0$ and concludes the proof of Lemma \ref{wxwx}.
\subsection{Coercivity of $H$} We complement the coercivity property of the operator $H^2$, established in the previous section, by the corresponding statement for the operator H, which follows from standard compactness argument. A complete proof is given in \cite{HR} with a slighlty different orthogonality condition but the proof is the same and therefore let to the reader.
\begin{lemma} [Coercivity of $H$]
\label{coercivityH}
Let $M \geq 1$ fixed. Then there exists $c(M) >0$ such that the following holds true. Let $u \in H^2_{rad}$ with
\bee
(u ,\Phi_M) = 0
\eee
then
\be
\int |\pa_{yy} u|^2 + \int \frac{|\pa_y u|^2}{y^2} + \int \frac{u^2}{y^4(1+|\log y|)^2} \leq c(M) \int |Hu|^2.
\ee
\end{lemma}
\section{Interpolation estimates}
In this section, we prove interpolation estimates for $\varepsilon$ in the bootstrap regimes which are used all along the proof of Proposition \ref{bootstrapregime}. We recall the norm $\mathcal E_1$, $\mathcal E_2$ and $\mathcal E_4$, introduced in \fref{defnorme}, together with their bootstrap bounds:
\bee
\mathcal E_1 &=& \int |\nabla \varepsilon|^2  \leq K \delta(b^*), \\
\mathcal E_2 &=& \int |H \varepsilon|^2 \leq Kb^2(t)|\log b(t)|^5,\\
\mathcal E_4 &=& \int |H^2 \varepsilon|^2 \leq K \frac{b^4(t)}{|\log b(t)|^2}.
\eee

\begin{lemma}[Interpolation estimates]
\label{lemmainterpolation}
There holds -with constants a priori depending on M-:
\bea
\label{estun}
  &&\int \frac{|H(\e)|^2}{y^4(1+|\log y|^2)} +  \int \frac{|\pa_y H(\e)|^2}{y^2(1+|\log y|^2)} +  \int \frac{|\pa_{yy} H(\e)|^2}{(1+|\log y|^2)}\\
 \nonumber &+& \int \frac{|\e|^2}{y^{4}(1+y^4)(1+|\log y|^2)} + \int \frac{|\pa^i_y\e|^2}{y^{8-2i}(1+|\log y|^2)} \lesssim \mathcal E_4, \ 1\leq i\leq 4,
\eea
\bea
\label{estunbis}
 \int \frac{|\e|^2}{y^{4}(1+|\log y|^2)} + \int \frac{|\pa^i_y\e|^2}{y^{4-2i}}\lesssim \mathcal E_2,\ \ 1\leq i\leq 2,
\eea
\be
\label{lossyboundwperp}
\int_{y \geq 1} \frac{1+|\log y|^C}{y^{8-2i}(1+|\log y|^2)}|\pa_y^i\e|^2 \lesssim b^4|\log b|^{C_1(C)}, \ \ 0\leq i\leq 2,
\ee
\be
\label{inteproloatedbound}
\int_{y \geq 1} \frac{1+|\log y|^C}{y^{6-2i}(1+|\log y|^2)}|\pa_y^i\e|^2 \lesssim b^3|\log b|^{C_1(C)}, \ \ 0\leq i\leq 2,
\ee
\be
\label{lossyboundwperpbis}
\int_{y \geq 1} \frac{1+|\log y|^C}{y^{4-2i}(1+|\log y|^2)}|\pa_y^i\e|^2 \lesssim b^2|\log b|^{C_1(C)}, \ \ 0\leq i\leq 1,
\ee
\be
\label{epsilonborne}
 \left \| \e (1+|\log y|^C) \right \|^2_{L^{\infty}_{y \geq 1}} \lesssim b^2|\log b|^{C_1(C)},
 \ee
\be
\label{epsilonborneavant1}
 \left \| \e \right \|^2_{L^{\infty}_{y \leq 1}} +  \left \| \pa_y \e \right \|^2_{L^{\infty}_{y \leq 1}}+  \left \|y \pa_{yy} \e \right \|^2_{L^{\infty}_{y \leq 1}} \lesssim \frac{b^4}{|\log b|^2},
 \ee
 \be
 \label{ypayepsilonborne}
 \left \|y \pa_y \e  \right \|^2_{L^{\infty}} \lesssim b^2|\log b|^5 ,
 \ee
\be
\label{normeinfty}
\left \| \frac{\varepsilon}{1+y} \right\|^2_{L^{\infty}} + \left\| \partial_y \varepsilon \right\|^2_{L^{\infty}} \lesssim b^3 |\log b|^C,
\ee
\be
\label{normeinfty2}
\left \| \frac{\varepsilon}{1+y^2} \right\|^2_{L^{\infty}} + \left\| \frac{\partial_y \varepsilon}{1+y} \right\|^2_{L^{\infty}} + \left\| \partial_{yy} \e \right\|^2_{L^{\infty}_{y \geq 1}} \lesssim b^4 |\log b|^C.
\ee

\end{lemma}
{\it Proof of the Lemma \ref{lemmainterpolation}:}
\fref{estun} and \fref{estunbis} are respectively direct consequences of the Lemma \ref{coercivityH2} and Lemma \ref{coercivityH} and  definition of the norms $\mathcal E_2$ and $\mathcal E_4$.  \\
To prove \fref{lossyboundwperp}, we split the integral at $y = B_0^{20}$.
\bee
&&\int_{y \geq 1} \frac{1+|\log y|^C}{y^{8-2i}(1+|\log y|^2)}|\pa_y^i\e|^2 \\
&=& \int_{1 \leq y \leq B_0^{20}} \frac{1+|\log y|^C}{y^{8-2i}(1+|\log y|^2)}|\pa_y^i\e|^2 + \int_{y \geq B_0^{20}} \frac{1+|\log y|^C}{y^{8-2i}(1+|\log y|^2)}|\pa_y^i\e|^2 \\
&\lesssim& |\log b|^{C+2} \mathcal E_4 \\
&+& \left \| \frac{1+|\log y|^C}{y^2}\right \|_{L^{\infty}_{y \geq B_0^{20}}} \left( \int_{y \geq 1} \frac{|\pa^i_y\e|^2}{y^{8-2i}(1+|\log y|^2)} \right)^\frac12 \left(  \int_{y \geq 1} \frac{|\pa^i_y\e|^2}{y^{4-2i}(1+|\log y|^2)}\right)^\frac12.
\eee
The bounds \fref{estun} and \fref{estunbis} concludes the proof. The bound \fref{inteproloatedbound} is a direct consequence of the last bound and the bootstrap bound for $\mathcal E_2$. Indeed:
\bee
&& \int_{y \geq 1} \frac{1+|\log y|^C}{y^{6-2i}(1+|\log y|^2)}|\pa_y^i\e|^2 \\
&\lesssim& \left( \int_{y \geq 1} \frac{1+|\log y|^{2C}}{y^{8-2i}(1+|\log y|^2)}|\pa_y^i\e|^2 \right) ^{\frac12}\left(  \int \frac{|\pa^i_y\e|^2}{y^{4-2i}(1+|\log y|^2)}\right) ^{\frac12}
\eee
and \fref{inteproloatedbound} follows. \par
The proof of \fref{lossyboundwperpbis} is the same of \fref{lossyboundwperp} using the energy $\mathcal E_1$ and $\mathcal E_2$.
\fref{epsilonborne} comes from \fref{normeinftysup1} and \fref{lossyboundwperpbis}. Indeed:
\\
\bee
 &&\left \| \e(1+|\log y|^C) \right \|^2_{L^{\infty}_{y \geq 1}}\\
& \lesssim& \int_{y \geq 1} \frac{|\pa_y \e|^2(1+|\log y|^C)}{y^2} + \int_{y \geq 1} \frac{\e^2(1+|\log y|^C)}{y^4} + \int_{\frac12 \leq y \leq 1} \frac{\e^2(1+|\log y|^C)}{y^2}\\
&\lesssim& b^2 |\log b|^C.
\eee
Let us prove \fref{epsilonborneavant1}. Let $a \in [1;2]$ such that:
\be
|\e(a)| \lesssim \int_1^2 \frac{|\e(y)|}{y^3} \lesssim \sqrt{\mathcal E_4}.
\ee
Then using Cauchy-Schwarz
\bee
\forall y \in [0;1], |\e(y)| \lesssim |\e(a)| + \left|\int_a^y \frac{|\pa_y \e(y)|}{y^3}\right| \lesssim  \sqrt{\mathcal E_4}.
\eee
In the same way, let $a \in [1;2]$ such that:
\be
|\pa_y\e(a)| \lesssim \int_1^2 \frac{|\pa_y\e(y)|}{y^3} \lesssim \sqrt{\mathcal E_4}
\ee
and
\bee
\forall y \in [0;1], |\pa_y \e(y)| \lesssim |\pa_y\e(a)| + \left|\int_a^y \frac{|\pa_{yy} \e(y)|}{y^3}\right| \lesssim  \sqrt{\mathcal E_4}.
\eee
Finally, let $a \in [1;2]$ such that:
\be
| \pa_{yy}\e(a)| \lesssim \int_1^2 \frac{|\pa_{yy}\e(y)|}{y^2} \lesssim \sqrt{\mathcal E_4}
\ee
and
\bee
\forall y \in [0;1], |y \pa_{yy} \e(y)| \lesssim a|\pa_{yy}\e(a)| + \left|\int_a^y\left( \frac{|\pa_{yy} \e(y)|}{y^3} +  \frac{|\pa^3_{y} \e(y)|}{y^2}\right)\right| \lesssim  \sqrt{\mathcal E_4}.
\eee
The bound \fref{ypayepsilonborne} is a direct consequence of the Lemma \ref{lemmetropcool}, and \ref{estunbis}.
\par
We now prove \fref{normeinfty} using \fref{normeinftysup1},  \fref{estunbis} and \fref{lossyboundwperp}:
\bee
&&\left \| \frac{\varepsilon}{y} \right\|^2_{L^{\infty}_{y \geq 1}} + \left\| \partial_y \varepsilon \right\|^2_{L^{\infty}_{y \geq 1}} \lesssim \int_{y \geq 1} \left( \frac{\e^2}{y^6} + \frac{|\pa_y \e|^2}{y^4} + \frac{|\pa_{yy} \e|^2}{y^2} \right)  + \int_{\frac12 \leq y \leq1} \left( \frac{|\e|^2}{y^4} + \frac{|\pa_y \e |^2}{y^2}\right)\\
&\lesssim&  \sum_{i=0}^{2} \left( \int_{y \geq 1} \frac{|\pa_y^i \e|^2}{y^{4-2i}(1+|\log y|^2)}\right)^{\frac12}\left( \int_{y \geq 1} \frac{|\pa_y^i \e|^2(1+|\log y|^2)}{y^{8-2i}} \right)^{\frac12} + \mathcal E_4 \lesssim b^3 |\log b|^C.
\eee
Similarily, we prove \fref{normeinfty2} :
\bee
&&\left \| \frac{\varepsilon}{y^2} \right\|^2_{L^{\infty}_{y \geq 1}} + \left \| \frac{\pa_y \varepsilon}{y} \right\|^2_{L^{\infty}_{y \geq 1}} + \left\| \partial_{yy} \varepsilon \right\|^2_{L^{\infty}_{y \geq 1}} \\
&\lesssim& \int_{y \geq 1} \left( \frac{\e^2}{y^8} + \frac{|\pa_y \e|^2}{y^6} + \frac{|\pa_{yy} \e|^2}{y^4} + \frac{|\pa^3_y \e|^2}{y^2} \right)  + \int_{\frac12 \leq y \leq1} \left( \frac{|\e|^2}{y^6} + \frac{|\pa_y \e |^2}{y^4}+ \frac{|\pa_{yy} \e |^2}{y^2}\right)\\
 &\lesssim& b^4 |\log b|^C.
\eee
\section{Localization of the profile}
In this appendix, we are going to give the important steps of the proof of the Proposition \ref{localisation} and \ref{localisationbis}. 
\\
\par
To begin, remark that the definition \fref{eqerreurloc} of $\tilde{\Psi}_b $ in the localization near $B_1$ gives two types of error. One is the result of the only localization. One is the effect of the time derivate. Indeed, we can rewrite $\tilde{\Psi}_b $ as following :
\be
\tilde{\Psi}_b  = \Psi_b^{(1)} + \tilde R
\ee
where
\be
\Psi_b^{(1)} = -b^2(\tilde T_1 + b \tilde T_2) - \Delta \tilde Q_b + b \Lambda \tilde Q_b - (\tilde Q_b)^3
\ee
and 
\be
\tilde R = b_s\left( 3b^2\tilde{T}_3 + b \frac{\pa \tilde T_1}{\pa b} + b^2 \frac{\pa \tilde T_2}{\pa b} + b^3 \frac{\pa \tilde T_3}{\pa b}\right).
\ee
We compute the action of localization which produces an error localized in $[B_1,2B_1]$ up to the term $(1 - \chi_{B_1})\Lambda Q$:
\bea
\label{defpsibloc}
 \nonumber \Psi^{(1)}_b &=& \chi_{B_1} \Psi_b + b(1 - \chi_{B_1})\Lambda Q +  b \Lambda \chi_{B_1}\alpha- 
\alpha \Delta \chi_{B_1} - 2 \partial_y \chi_{B_1} \partial_y\alpha\\
& + & (Q+\chi_{B_1}\alpha)^3-Q^3-\chi_{B_1}((Q+\alpha)^3-Q^3).
 \eea
 We estimate from the rough bounds of $\left(T_i \right)_{1 \leq i \leq 3}$ and the choice of $B_1$: 
 $$\forall y\leq 2B_1, \  \ |\alpha(y)|\lesssim by\left(\frac{|\log y|}{y}+\frac{by}{|\log b|}\right)\lesssim b|\log y|$$ and thus:
 $$\left|b(1 - \chi_{B_1})\Lambda Q +  b \Lambda \chi_{B_1}\alpha- 
\alpha \Delta \chi_{B_1} - 2 \partial_y \chi_{B_1} \partial_y\alpha\right|\lesssim \frac{b}{y^2}{\bf 1}_{y\geq B_1}+b^2\log y {\bf 1}_{B_1\leq y\leq 2B_1},
$$
 \bee
 \left| (Q+\chi_{B_1}\alpha)^3-Q^3-\chi_{B_1}((Q+\alpha)^3-Q^3)\right| & \lesssim & \frac{|\alpha(y)|}{y^3}|{\bf 1}_{B_1\leq y\leq 2B_1}.\\
 & \lesssim & \frac{b\log y}{y^2}{\bf 1}_{B_1\leq y\leq 2B_1}
 \eee
Hence, $\Psi_b^{(1)}$ verifies the bounds \fref{controleh2erreurtilde}, \fref{newbornestilde}, \fref{controleh4erreurtilde}. For the control of time derivates, we have to use that :
\be
\label{fgtr}
\frac{\pa c_b}{\pa b} = O \left( \frac{1}{b|\log b|^2}\right), \ \ \frac{\pa d_b}{\pa b} = O \left( \frac{1}{b^2|\log b|^2}\right),
\ee
which are consequences of their definition, and that :
\be
\label{frtg}
\frac{\pa \Sigma_b}{\pa b} = O \left( \frac{1}{b|\log b|} {\bf 1}_{y \leq \frac{B_0}{2}} + \frac{1}{y^2b^2|\log b|^2} {\bf 1}_{\frac{B_0}{2} \leq y \leq 6B_0} \right).
\ee
Using \fref{fgtr}, \fref{frtg} together the explicit formula of $\left(T_i \right)_{1 \leq i \leq 3}$ yield \fref{controleh2erreurtilde}, \fref{newbornestilde}, \fref{controleh4erreurtilde} without difficulty. Only an estimate is more delicate, and requests more cancellation. Indeed, we must use that 
\bee
H^2\left( \frac{\pa \tilde T_2}{\pa b}\right) = H \left( \frac{\pa \Sigma_2}{\pa b}\right) = H\left( \frac{\pa \Sigma_b}{\pa b}\right),
\eee
and similarly 
\bee
H^2\left( \frac{\pa \tilde T_3}{\pa b}\right) = H \left( \frac{\pa \Sigma_3}{\pa b}\right) = \Lambda \left( \frac{\pa \Sigma_b}{\pa b}\right).
\eee
The proof of \fref{fluxcomputationonebis} is left to the reader.
For the proof of Proposition \ref{localisationbis}, we just remark that, by definition : $$\zeta_b=(\chi_{B_1}-\chi_{B_0})(bT_1+b^2T_2+b^3T_3)$$
This proof is afterwards the same as the previous one.
\end{appendix}

\end{document}